\documentclass[12pt]{article}
\usepackage{amssymb}
\usepackage{amsmath}
\usepackage{graphics}
\usepackage{color}
\renewcommand{\baselinestretch}{1.2}

\newcounter{theorem}[section]
\newcounter{remark}[section]
\newcounter{corollary}[section]
\newcounter{lemma}[section]
\newcounter{proposition}[section]
\newcounter{definition}[section]
\newcounter{assumption}
\newcounter{notice}[section]

\def\vecc#1{\mathord{\buildrel{\lower3pt\hbox{$\scriptscriptstyle\rightharpoonup$}}\over #1}}

\begin{document}
\title{\bf HFVS: An Arbitrary High Order Flux Vector Splitting Method
 \footnote{Supported by
NSFC(Grant No. 11101047) .} }
\author{{ Yibing Chen\ , Song Jiang \ and\ \  Na Liu}\\
{\small    Institute of Applied Physics and Computational Mathematics,}\\[-2mm]
{\small    P.O. Box 8009, Beijing 100088, P.R. China}\\[-1mm]
{\small    E-mail: chen\_{}yibing@iapcm.ac.cn,\ \
jiang@iapcm.ac.cn,\ \ liu\_{}na@iapcm.ac.cn} }
\date{ }

\maketitle
\begin{abstract}
In this paper, a new scheme of arbitrary high order accuracy in both
space and time is proposed to solve hyperbolic conservative laws.
Based on the idea of flux vector splitting(FVS) scheme, we split all
the space and time derivatives in the Taylor expansion of the
numerical flux into two parts: one part with positive eigenvalues,
another part with negative eigenvalues. According to a Lax-Wendroff
procedure, all the time derivatives are then replaced by space
derivatives. And the space derivatives is calculated by  WENO
reconstruction polynomial. One of the most important advantages of
this new scheme is easy to implement.In addition, it should be
pointed out, the procedure of calculating the space and time
derivatives in numerical flux can be used as a building block to
extend the current first order schemes to very high order accuracy
in both space and time. Numerous numerical tests for linear
 and nonlinear hyperbolic conservative laws demonstrate
that new scheme is robust and can be high order accuracy in both
space and time.

\end{abstract}


\section{Introduction}
In the recent years, numerical schemes of high order accuracy for
hyperbolic conservative laws have attracted much attention.
In the 1980s, thanks to the efforts of Harten \cite{harten84}, Van
Leer \cite{vanLeer77} and Roe \cite{Roe81} et al. (only a partial list
is given here), schemes of second order accuracy tend to maturity.
Later, ENO/WENO schemes show a possibility to construct the schemes of
arbitrary high order accuracy in space \cite{Harten87,weno96}.
And since the end of the 1990s, designing schemes of high order accuracy (greater
than second order accuracy) becomes one of the most important issues in CFD.
In practical simulations, ENO/WENO schemes usually couple with the
multi-stage Runge-Kutta method for time evolution.
In implementation , the third order TVD Runge-Kutta method is often
used \cite{bs00}. This is because higher than third order
Runge-Kutta methods will become complicate. For example, fourth
order TVD Runge-Kutta method needs to save all the variables in
intermediate stages, which will definitely introduce more
computational complexity. What's worse, higher than fourth order
Runge-Kutta method will meet the so-called Butcher
barrier\cite{Butcher}, which the number of stages will greater than
the order of accuracy (e.g., fifth order Runge-Kutta method need six
stages). But third order Runge-Kutta method will also lead to
accuracy barrier, for example, the fifth order WENO scheme couples
with the third order TVD Runge-Kutta method can only achieve third
order accuracy \cite{Toro02}.
To obtain ideal convergence rate, one has to reduce the CFL
condition number, which will increase the computational costs.

At the first decade of this century, many researchers made great
efforts to avoid the above drawback. The first milestone may belong
to the ADER (Arbitrary DERivative in space and time) scheme, which
is a Godunov approach based scheme of arbitrary high order accuracy
in both space and time \cite{Toro02}. 
The idea of ADER scheme can go back to the GRP scheme \cite{grp84},
which solves the generalized Riemann problem with piecewise linear
initial values instead of the conventional Riemann problem with piecewise constant initial values.
The GRP scheme is sophisticated, but may become quite complicated in the case of higher than third order accuracy.
By introducing a linearization technique, the ADER scheme transfers a
generalized Riemann problem into a conventional Riemann problem
coupled with a series of linear hyperbolic problems. This transfer process
much simplifies the computational process of the original GRP scheme
and makes it possible to construct schemes of arbitrary high order
accuracy both in space and time.

The ADER scheme can be divided into two camps, the one is
state-expansion \cite{Toro02}, another one is flux-expansion \cite{Toro05}.
The former can not be integrated explicitly in time, so a Gaussian quadrature
is necessary for high order accuracy.
On the other hand, the latter can be integrated exactly in time, the
time derivatives of the numerical flux, however, have to be given
explicitly, and it could be very tedious when the formulation of the
numerical flux is complex. If one chooses a nonlinear Riemann
problem solver such as HLLC, it is not trivial to deduce the
formulation of high order derivative terms.
\textcolor[rgb]{1.00,0.00,0.00}{That is why the flux expansion
version ADER is more difficult to implement than the state expansion
one.}

Besides ADER scheme, HGKS (high-order accurate gas-kinetic scheme)
is another important scheme of high order accuracy
\cite{hgks10,hgks14}, which is a gas kinetic scheme.
HGKS is based on the Boltzmann equation, then the time derivatives
in the Taylor expansion of the particle distribution function can be
replaced by space derivatives directly.
Comparing to the ADER scheme, the numerical flux in HGKS can be obtained without a linearization process.
HGKS can be also integrated exactly in time and made to be of arbitrary high order
accuracy both in space and time. We should point out, however, that the formulation of HGKS may also
become quite complex in the high order accuracy case. Although the recursive techniques simplify the
procedure of deduction in the high order HGKS, the computational costs can not be reduced. In fact,
the computational costs in CPU time for HGKS may greater than those for WENO with Godunov
Riemann problem solvers, see \cite{hgks13}.
%
Another drawback for HGKS is that there are some spurious velocity
and pressure oscillations near the contact discontinuous when HGKS
is used to solve the Euler equations. This phenomenon is found and
carefully analyzed by the authors in \cite{cj09}. A remedy for the
GKS coupled with the Runge-Kutta method is also provided in their
articles \cite{cj09,cj11}. However, how to construct an
oscillation-free HGKS (greater than second order accuracy) is still
unsolved.

In recent fifteen years, the ADER scheme and HGKS have provided a
new approach to construct schemes of high order accuracy. In this
article, we present a new simple version to construct a high order
scheme and propose thus a numerical scheme of arbitrary high order
accuracy in both space and time.
Different from the ADER scheme and HGKS, our new scheme is based on
the flux vector splitting method, and is much easier to implement
and robust. We call the new scheme HFVS (arbitrary high order flux
vector splitting scheme). We shall give a number of numerical tests
to demonstrate the efficiency and accuracy of HFVS.

This paper is organized as follows. In Section 2, we present the
framework of numerical schemes of arbitrary high order accuracy in both space
and time, while in Section 3 we give the construction of our
new scheme HFVS in details.
In Section 4 we test HFVS in accuracy by a number of numerical examples of linear and nonlinear
conservative laws in 1D and 2D, and compare the numerical results computed by HFVS and those by WENO in both
accuracy and computational costs.
Finally, we draw the conclusions.

\section{Framework of numerical schemes of high order accuracy in both space and time}

As mentioned in the previous section, the ADER scheme is a Godunov
scheme based approach, while HGKS is a Boltzmann equation based
approach. In our opinion, the flux-expansion version of the ADER
scheme and HGKS can be classified in a unified framework of high
order schemes in both space and time.
\textcolor[rgb]{1.00,0.00,0.00}{In fact, Toro's original
article\cite{Toro05} implies the possibility. Here we want to make
the framework more clearly.}

Let us start with a one-dimensional hyperbolic system of conservation laws in the form:
\begin{equation} \label{eq:1}
\frac{\partial\vecc{W}}{\partial t} + \frac{\partial \vecc{F}(\vecc{W})}{\partial x} = 0.
\end{equation}
To simplify the presentation, we consider here the 1D Euler equations of compressible
flows in this section, although both the ADER scheme and
HGKS are constructed for general hyperbolic conservation laws. For the the 1D Euler equations,
$\vec{W}$ and $F$ in (\ref{eq:1}) are given by
\begin{equation} \label{1deuler}
\begin{array}{ll}
\vecc{W}  = \big( {\rho ,\rho U,\rho E} \big)^T , \\[2mm]
\vecc{F}(\vecc{W}) = \big( {\rho U,\rho U^2  + P,\rho EU + PU} \big)^T, \end{array}
\end{equation}
where, $\rho,\, U,\, E,\, P$ are the density, fluid velocity, specific total
energy and pressure, respectively. The system can be closed by
adding the equations of state (perfect gases):
$$ P = (\gamma  - 1)\rho\Big( E - \frac{1}{2}U^2\Big). $$

 Integrating (\ref{eq:1}) over a space-time computational cell
$[{x_{j - \frac{1}{2}}},{x_{j + \frac{1}{2}}}] \times [{t^n},{t^{n + 1}}]$, we have
\begin{equation} \label{eq:2}
\int_{{t^n}}^{{t^{n + 1}}} {\int_{{x_{j - \frac{1}{2}}}}^{{x_{j + \frac{1}{2}}}} {\left( {\frac{{\partial
\mathord{\buildrel{\lower3pt\hbox{$\scriptscriptstyle\rightharpoonup$}}
\over W} }}{{\partial t}} + \frac{{\partial
\mathord{\buildrel{\lower3pt\hbox{$\scriptscriptstyle\rightharpoonup$}}
\over F}
(\mathord{\buildrel{\lower3pt\hbox{$\scriptscriptstyle\rightharpoonup$}}
\over W} )}}{{\partial x}}} \right)} } dxdt = 0.
\end{equation}

Denote $\Delta x = [{x_{j - \frac{1}{2}}},{x_{j + \frac{1}{2}}}]$
and $\Delta t = [{t^n},{t^{n + 1}}]$, we approximate (\ref{eq:2}) by the
classical finite volume method (FVM)
\begin{equation} \label{eq:3}
\mathord{\buildrel{\lower3pt\hbox{$\scriptscriptstyle\rightharpoonup$}}
\over W} _j^{n + 1} =
\mathord{\buildrel{\lower3pt\hbox{$\scriptscriptstyle\rightharpoonup$}}
\over W} _j^n - \frac{1}{{\Delta x}}\int_{{t^n}}^{{t^{n + 1}}}
{\left( {{{\mathord{\buildrel{\lower3pt\hbox{$\scriptscriptstyle\rightharpoonup$}}
\over F} }_{j + \frac{1}{2}}} -
{{\mathord{\buildrel{\lower3pt\hbox{$\scriptscriptstyle\rightharpoonup$}}
\over F} }_{j - \frac{1}{2}}}} \right)} dt ,
\end{equation}
where $\mathord{\buildrel{\lower3pt\hbox{$\scriptscriptstyle\rightharpoonup$}}\over W}_j^n$ is the cell average of
$\mathord{\buildrel{\lower3pt\hbox{$\scriptscriptstyle\rightharpoonup$}}
\over W} (x,t^n )$ ($x \in [x_{j - \frac{1}{2}},x_{j+\frac{1}{2}}]$) and
$\mathord{\buildrel{\lower3pt\hbox{$\scriptscriptstyle\rightharpoonup$}}\over F}_{j+\frac{1}{2}}$
is the numerical flux. To make the scheme (\ref{eq:3}) high order accuracy, the key point is to suitably construct the numerical
flux which should be in high order accuracy.

The procedure to construct a numerical scheme with high order
accuracy in both space and time can be divided into three steps.

\textbf{Step I.} Construct $\mathord{\buildrel{\lower3pt\hbox{$\scriptscriptstyle\rightharpoonup$}}\over F}_{j+\frac{1}{2}}(0^+)$
by the Taylor expansion of the numerical flux as follows.
\begin{equation} \label{eq:4}
{\mathord{\buildrel{\lower3pt\hbox{$\scriptscriptstyle\rightharpoonup$}}
\over F}_{j + \frac{1}{2}}}(\tau ) =
{\mathord{\buildrel{\lower3pt\hbox{$\scriptscriptstyle\rightharpoonup$}}
\over F}_{j + \frac{1}{2}}}({0^ + }) + \sum\limits_{k = 1}^N
{\frac{{{\partial ^k}}}{{\partial {t^k}}}}
{\mathord{\buildrel{\lower3pt\hbox{$\scriptscriptstyle\rightharpoonup$}}
\over F}_{j + \frac{1}{2}}}({0^ + })\frac{{{\tau ^k}}}{{k!}}.
\end{equation}

\textbf{Step II.} Replace the time derivatives by spatial derivatives.

\textbf{Step III.} Calculate the space derivatives by certain reconstruction techniques.

In the above procedure, a WENO reconstruction technique is used in
Step III for the ADER scheme and HGKS, while very different techniques are employed in Steps I and II, and described
in the following.

\subsection{The ADER scheme}

In the ADER scheme, an approximate Riemann solver such as HLLC is used in Step II:
\[
\mathord{\buildrel{\lower3pt\hbox{$\scriptscriptstyle\rightharpoonup$}}
\over F} _{j + \frac{1}{2}} (0^ +  ) = \left\{
{\begin{array}{*{20}c}
   {\mathord{\buildrel{\lower3pt\hbox{$\scriptscriptstyle\rightharpoonup$}}
\over F} _j ,\qquad \quad \quad \quad \quad if\;0 \le \;S_j },  \\
   {\mathord{\buildrel{\lower3pt\hbox{$\scriptscriptstyle\rightharpoonup$}}
\over F} _j  + S_j \left(
{\mathord{\buildrel{\lower3pt\hbox{$\scriptscriptstyle\rightharpoonup$}}
\over W} _j^*  -
\mathord{\buildrel{\lower3pt\hbox{$\scriptscriptstyle\rightharpoonup$}}
\over W} _j } \right),\;\;if\;S_j  \le 0 \le \;S^* },  \\
   {\mathord{\buildrel{\lower3pt\hbox{$\scriptscriptstyle\rightharpoonup$}}
\over F} _{j + 1}  + S_{j + 1} \left(
{\mathord{\buildrel{\lower3pt\hbox{$\scriptscriptstyle\rightharpoonup$}}
\over W} _j^*  -
\mathord{\buildrel{\lower3pt\hbox{$\scriptscriptstyle\rightharpoonup$}}
\over W} _j } \right),\;\; if\;S^*  \le 0 \le \;S_{j + 1} },  \\
   {\mathord{\buildrel{\lower3pt\hbox{$\scriptscriptstyle\rightharpoonup$}}
\over F} _{j + 1} ,\qquad \qquad \quad \quad if\;0\; \ge S_j },  \\
\end{array}} \right. \]
where
\[ \mathord{\buildrel{\lower3pt\hbox{$\scriptscriptstyle\rightharpoonup$}}
\over W} _k^*  = \rho _k \left( {\frac{{S_k  - U_k }}{{S_k  - U^*
}}} \right)\left[ {\begin{array}{*{20}c}    1  \\
   {S^* }  \\
   {E_k  + \left( {S^*  - U_k } \right)\left( {S^*  + \frac{{P_k }}{{\rho _k \left( {S_k  - U_k } \right)}}} \right)}  \\
\end{array}} \right],k = j\;or\;j + 1  \]
The wave speed of $S_j$, $S_j+1$ and $S^*$ can have many choices, see \cite{hllc}.

Then, the time derivatives of numerical flux in (\ref{eq:4}) can be given analytically, for example, as follows.
\begin{eqnarray*} &&
 \frac{{\partial \mathord{\buildrel{\lower3pt\hbox{$\scriptscriptstyle\rightharpoonup$}}
\over F} }}{{\partial t}} = \Big[ {\left( {\rho U} \right)_t ,\left( {\rho U} \right)_t U + \left( {\rho U} \right)U_t
+ P_t ,U_t \left( {\rho E + P} \right) + U\left( {\rho E + P} \right)_t } \Big]^T,  \\
&& \frac{{\partial ^2 \mathord{\buildrel{\lower3pt\hbox{$\scriptscriptstyle\rightharpoonup$}}
\over F} }}{{\partial t^2 }} = \Big[ (\rho U)_{tt} ,(\rho U)_{tt} U + 2(\rho U)_t U_t
 + (\rho U)U_{tt}  + P_{tt} ,U_{tt}(\rho E + P)  \\
 && \qquad + 2U_t(\rho E + P)_t  + U(\rho E + P)_{tt} \Big]^T.
 \end{eqnarray*}

According to the construction procedure of the Lax-Wendroff scheme, the time derivatives in the
above formulations can be replace by the spatial derivatives.

\textbf{Remark 1.} The procedure can be carried out for any given
order. However, it is easy to see that these formulations may become
very complicate in the case of very high order accuracy , which will
increase the complexity in implementation.

\subsection{HGKS}
Different from the ADER scheme, HGKS is based on the Boltzmann-type equation, such as the BGK model:
\[ f_t  + uf_x  = \frac{{g - f}}{\tau }, \]
where $f$ is the particle velocity distribution function and $g$ is
corresponding equilibrium distribution. Both $f$ and $g$ are
functions of space $x$, time $t$, the particle velocity $u$, and the
internal variables $\vec{\xi}=(\xi_1,\xi_2,\cdots,\xi_K)$, which has $K$ degrees of freedom.

The relationship between the numerical flux and the particle distribution function is
\[
\mathord{\buildrel{\lower3pt\hbox{$\scriptscriptstyle\rightharpoonup$}}
\over F} _{j + \frac{1}{2}} \left( \tau  \right) = \int {u\Psi f_{j
+ \frac{1}{2}} \left( \tau  \right)} d\Xi.
\]
Here $\Psi$ is the moment vector
\[
\Psi=(1,u,\frac12(u^2+\overrightarrow{\xi}^2))^T,
\]
and the internal variable $\vec{\xi}^2$ is equal to $\vec{\xi}^2=\xi_1^2+\xi_2^2+\cdots +\xi_K^2$.

Hence, the time derivatives of numerical flux can be calculated by
\[
\frac{{\partial ^k \mathord{\buildrel{\lower3pt\hbox{$\scriptscriptstyle\rightharpoonup$}}
\over F} _{j + \frac{1}{2}} }}{{\partial t^k }} = \int {u\Psi
\frac{{\partial ^k f_{j + \frac{1}{2}} }}{{\partial t^k }}} d\Xi . \]

In view of the BGK model, the time derivatives can be evaluated directly, for example,
\begin{eqnarray*} &&
f _{j + \frac{1}{2}} (0^ +  )= (1-\tau(a^l u+A^l))H(u)g^l+ (1-\tau(a^r u+A^r))(1-H(u))g^r,  \\[1mm]
&& \frac{\partial}{\partial t}f_{j+\frac12}(0^+)=\frac1\tau g_0-\tau(a^2+b)u^2g_0-\tau(A^2+B')g_0-2\tau(C+aA)ug_0 \\
&& \quad +[-\frac1\tau+A^l+\tau(C^l+a^lA^l)u+\tau((a^l)^2+b^l)u^2]H(u)g^l \\
&& \quad +[-\frac1\tau+A^r+\tau(C^r+a^rA^r)u+\tau((a^r)^2+b^r)u^2](1-H(u))g^r, \\[1mm]
&& \frac{\partial^2}{\partial t^2}f_{j+\frac12}(0^+)=-\frac{1}{\tau^2} g_0-\frac{1}{\tau}aug_0+\frac{1}{\tau}Ag_0+(a^2+b)u^2g_0+(A^2+B')g_0+2(C+aA)ug_0\\
&& \quad +[\frac{1}{\tau^2}-\frac{1}{\tau}(a^lu+A^l)+\frac{2}{\tau}a^lu-2(C^l+a^lA^l)u-((a^l)^2+b^l)u^2]H(u)g^l \\
&& \quad +[\frac{1}{\tau^2}-\frac{1}{\tau}(a^ru+A^r)+\frac{2}{\tau}a^ru-2(C^r+a^rA^r)u-((a^r)^2+b^r)u^2](1-H(u))g^r.
\end{eqnarray*}
Here $g_0$ is the Maxwellian distribution at the cell interface
$x_{j+\frac12}$, $g^l$ and $g^r$ are the left and right limits of
the Maxwellian distribution at the cell interface, repsectively.
$a,b,A,\cdots$ are related to the space and time derivatives of the
the Maxwellian distribution, see \cite{hgks14} for the detailed
definition.

\textbf{Remark 2.} Theoretically, HGKS can achieve arbitrary high
order accuracy. However, due to its complexity in the expression,
only third and fourth order schemes are often used in
simulations\cite{hgks10,hgks14}.

In the following section, we want to present a simpler high order scheme, which is different from the ADER scheme and HGKS.

\section{HFVS: Arbitrary high order flux vector splitting scheme}

%

From the last section, we can see that the complexity of the flux expansion version of the ADER scheme
and HGKS actually comes from the procedure of evaluating the time derivatives of the numerical fluxes.
So, if one can find a numerical flux, the time derivatives of which can be calculated very easily,
then the complexity of the associated numerical scheme can be largely reduced. This motivates us to construct our scheme.

As is well-known, in comparison with the traditional Godunov scheme,
flux vector splitting (FVS) schemes are simpler in construction and
coding. Thus, we try to use the FVS approach, instead of the Godunov
approach used in the ADER scheme. Our proposed scheme is constructed
in the same framework as in Section 2, and also consists of three
phases: constructing spatial derivatives, leading terms and time
derivatives.

\subsection{Spatial derivatives}
In order to diminish possible spurious oscillations, we should use
conservative reconstruction techniques such as ENO/WENO to evaluate
spatial derivatives.

Suppose that a smooth function over a cell should be reconstructed,
which is of the form:
\begin{equation} \label{eq:10}
\mathord{\buildrel{\lower3pt\hbox{$\scriptscriptstyle\rightharpoonup$}}
\over W} (x) =
\mathord{\buildrel{\lower3pt\hbox{$\scriptscriptstyle\rightharpoonup$}}
\over W}_j  + \sum\limits_{k = 1}^N {\frac{{\partial^k
\mathord{\buildrel{\lower3pt\hbox{$\scriptscriptstyle\rightharpoonup$}}
\over W} }}{{\partial x^k }}\phi_k (x)},\qquad x \in [{x_{j
-\frac{1}{2}},x_{j+\frac{1}{2}} }],
\end{equation}
where
$\mathord{\buildrel{\lower3pt\hbox{$\scriptscriptstyle\rightharpoonup$}}\over
W}_j $ is the cell known averaged value of
$\mathord{\buildrel{\lower3pt\hbox{$\scriptscriptstyle\rightharpoonup$}}\over
W}$ over cell $j$, $\partial^k
\mathord{\buildrel{\lower3pt\hbox{$\scriptscriptstyle\rightharpoonup$}}\over
W} /\partial x^k$ ($k = 1,\cdots, N$) are the unknowns to be
determined, and
\begin{eqnarray*} &&
\varphi _1  = \Big( {\frac{{x - x_j }}{{\Delta x}}} \Big),\quad
\varphi_2 = \frac{1}{2}\Big[ {\Big( {\frac{{x - x_j }}{{\Delta
x}}}\Big)^2  - \frac{1}{{12}}} \Big], \quad
\varphi_3  = \frac{1}{6}\Big({\frac{{x - x_j }}{{\Delta x}}} \Big)^3 , \\[2mm]
&& \varphi_4  = \frac{1}{{24}}\Big[ {\Big( {\frac{x - x_j}{\Delta
x}}\Big)^4  - \frac{1}{{80}}} \Big],\quad \cdots
\end{eqnarray*}

In general, we can always obtain two values on cell interfaces by
the ENO/WENO reconstruction techniques, i.e.,
$\mathord{\buildrel{\lower3pt\hbox{$\scriptscriptstyle\rightharpoonup$}}
\over W}_{j\pm\frac{1}{2}}$. Thus, we immediately obtain two linear
\textbf{algebraic} equations:
\begin{equation} \label{eq:11}
\left\{ {\begin{array}{*{20}c}
   {\mathord{\buildrel{\lower3pt\hbox{$\scriptscriptstyle\rightharpoonup$}}
\over W} (x_{j + \frac{1}{2}} ) =
\mathord{\buildrel{\lower3pt\hbox{$\scriptscriptstyle\rightharpoonup$}}
\over W} _{j + \frac{1}{2}} },  \\
   {\mathord{\buildrel{\lower3pt\hbox{$\scriptscriptstyle\rightharpoonup$}}
\over W} (x_{j - \frac{1}{2}} ) =
\mathord{\buildrel{\lower3pt\hbox{$\scriptscriptstyle\rightharpoonup$}}
\over W} _{j - \frac{1}{2}} }.  \\
\end{array}} \right.
\end{equation}

In this case , these formulations also implies
\begin{equation} \label{eq:12}
\frac{1}{{\Delta x}}\int_{I_j }
{\mathord{\buildrel{\lower3pt\hbox{$\scriptscriptstyle\rightharpoonup$}}
\over W} (x)dx}  =
\mathord{\buildrel{\lower3pt\hbox{$\scriptscriptstyle\rightharpoonup$}}
\over W}_j.
\end{equation}
 When the dimension $N=1$, the above system is uniquely
solvable.

For $N>1$, supplementary information should be added in order to
have a solution of the system (\ref{eq:11}). A simple way of giving
the supplementary information is to employ the cell averaged values.
For example, for $N=2$ we use
\begin{equation} \label{eq:12}
\frac{1}{{\Delta x}}\int_{I_j }
{\mathord{\buildrel{\lower3pt\hbox{$\scriptscriptstyle\rightharpoonup$}}
\over W} (x)dx}  =
\mathord{\buildrel{\lower3pt\hbox{$\scriptscriptstyle\rightharpoonup$}}
\over W}_j.
\end{equation}
Coupling (\ref{eq:12}) with (\ref{eq:11}), we obtain a system of
three independent linear equations for the three unknowns, which is
uniquely solvable.

Similarly, for $N=4$ we use
\[
\frac{1}{{\Delta x}}\int_{I_l }
{\mathord{\buildrel{\lower3pt\hbox{$\scriptscriptstyle\rightharpoonup$}}
\over W} (x)dx}  =
\mathord{\buildrel{\lower3pt\hbox{$\scriptscriptstyle\rightharpoonup$}}\over
W}_l, \qquad l = j-1,j,j + 1 ,
\]
instead of (\ref{eq:12}). Obviously, this choice is suitable for the
case of arbitrary high order accuracy.

In fact, the spatial derivatives can be given explicitly, see below.

\textbf{Second order.}
 \begin{equation} \label{eq:21}
\left( {\frac{{\partial
\mathord{\buildrel{\lower3pt\hbox{$\scriptscriptstyle\rightharpoonup$}}
\over W} }}{{\partial x}}} \right)_j  =
\mathord{\buildrel{\lower3pt\hbox{$\scriptscriptstyle\rightharpoonup$}}
\over W} _{j + \frac{1}{2}}  -
\mathord{\buildrel{\lower3pt\hbox{$\scriptscriptstyle\rightharpoonup$}}
\over W} _{j - \frac{1}{2}}
\end{equation}

\textbf{Third order.}
\begin{equation} \label{eq:22}
\begin{array}{l}
 {\left( {\frac{{\partial \mathord{\buildrel{\lower3pt\hbox{$\scriptscriptstyle\rightharpoonup$}}
\over W} }}{{\partial x}}} \right)_j} = {{\mathord{\buildrel{\lower3pt\hbox{$\scriptscriptstyle\rightharpoonup$}}
\over W} }_{j + \frac{1}{2}}} - {{\mathord{\buildrel{\lower3pt\hbox{$\scriptscriptstyle\rightharpoonup$}}
\over W} }_{j - 1}} \\
 {\left( {\frac{{{\partial ^2}\mathord{\buildrel{\lower3pt\hbox{$\scriptscriptstyle\rightharpoonup$}}
\over W} }}{{\partial {x^2}}}} \right)_j} = 3( - 2{{\mathord{\buildrel{\lower3pt\hbox{$\scriptscriptstyle\rightharpoonup$}}
\over W} }_{j + 1}} + {{\mathord{\buildrel{\lower3pt\hbox{$\scriptscriptstyle\rightharpoonup$}}
\over W} }_{j + \frac{1}{2}}} + {{\mathord{\buildrel{\lower3pt\hbox{$\scriptscriptstyle\rightharpoonup$}}
\over W} }_{j - \frac{1}{2}}}) \\
 \end{array}
\end{equation}

\textbf{Fifth order.}
\begin{equation} \label{eq:23}
\begin{array}{l}
 {\left( {\frac{{\partial \mathord{\buildrel{\lower3pt\hbox{$\scriptscriptstyle\rightharpoonup$}}
\over W} }}{{\partial x}}} \right)_j} = \frac{1}{8}({{\mathord{\buildrel{\lower3pt\hbox{$\scriptscriptstyle\rightharpoonup$}}
\over W} }_{j - 1}} - {{\mathord{\buildrel{\lower3pt\hbox{$\scriptscriptstyle\rightharpoonup$}}
\over W} }_j} - 10{{\mathord{\buildrel{\lower3pt\hbox{$\scriptscriptstyle\rightharpoonup$}}
\over W} }_{j - \frac{1}{2}}} + 10{{\mathord{\buildrel{\lower3pt\hbox{$\scriptscriptstyle\rightharpoonup$}}
\over W} }_{j + \frac{1}{2}}}) \\
 {\left( {\frac{{{\partial ^2}\mathord{\buildrel{\lower3pt\hbox{$\scriptscriptstyle\rightharpoonup$}}
\over W} }}{{\partial {x^2}}}} \right)_j} = \frac{1}{8}( - {{\mathord{\buildrel{\lower3pt\hbox{$\scriptscriptstyle\rightharpoonup$}}
\over W} }_{j - 1}} - 50{{\mathord{\buildrel{\lower3pt\hbox{$\scriptscriptstyle\rightharpoonup$}}
\over W} }_j} - {{\mathord{\buildrel{\lower3pt\hbox{$\scriptscriptstyle\rightharpoonup$}}
\over W} }_{j + 1}} + 30{{\mathord{\buildrel{\lower3pt\hbox{$\scriptscriptstyle\rightharpoonup$}}
\over W} }_{j + \frac{1}{2}}} + 30{{\mathord{\buildrel{\lower3pt\hbox{$\scriptscriptstyle\rightharpoonup$}}
\over W} }_{j - \frac{1}{2}}}) \\
 {\left( {\frac{{{\partial ^3}\mathord{\buildrel{\lower3pt\hbox{$\scriptscriptstyle\rightharpoonup$}}
\over W} }}{{\partial {x^3}}}} \right)_j} = \frac{1}{2}( - {{\mathord{\buildrel{\lower3pt\hbox{$\scriptscriptstyle\rightharpoonup$}}
\over W} }_{j - 1}} + {{\mathord{\buildrel{\lower3pt\hbox{$\scriptscriptstyle\rightharpoonup$}}
\over W} }_{j + 1}} - 2{{\mathord{\buildrel{\lower3pt\hbox{$\scriptscriptstyle\rightharpoonup$}}
\over W} }_{j + \frac{1}{2}}} + 2{{\mathord{\buildrel{\lower3pt\hbox{$\scriptscriptstyle\rightharpoonup$}}
\over W} }_{j - \frac{1}{2}}}) \\
 {\left( {\frac{{{\partial ^4}\mathord{\buildrel{\lower3pt\hbox{$\scriptscriptstyle\rightharpoonup$}}
\over W} }}{{\partial {x^4}}}} \right)_j} = \frac{5}{{12}}({{\mathord{\buildrel{\lower3pt\hbox{$\scriptscriptstyle\rightharpoonup$}}
\over W} }_{j - 1}} + 10{{\mathord{\buildrel{\lower3pt\hbox{$\scriptscriptstyle\rightharpoonup$}}
\over W} }_j} + {{\mathord{\buildrel{\lower3pt\hbox{$\scriptscriptstyle\rightharpoonup$}}
\over W} }_{j + 1}} - 6{{\mathord{\buildrel{\lower3pt\hbox{$\scriptscriptstyle\rightharpoonup$}}
\over W} }_{j + \frac{1}{2}}} - 6{{\mathord{\buildrel{\lower3pt\hbox{$\scriptscriptstyle\rightharpoonup$}}
\over W} }_{j - \frac{1}{2}}}) \\
 \end{array}
 \end{equation}

\subsection{Leading terms }
 Following the idea of FVS, the flux in a computational cell can be decomposed into two parts:
\[  \mathord{\buildrel{\lower3pt\hbox{$\scriptscriptstyle\rightharpoonup$}}
\over F}_j  = \mathord{\buildrel{\lower3pt\hbox{$\scriptscriptstyle\rightharpoonup$}}\over F}_j^+
 +\mathord{\buildrel{\lower3pt\hbox{$\scriptscriptstyle\rightharpoonup$}}\over F}_{j + 1}^- .  \]
Thus the numerical flux on a cell edge consists of two parts:
\begin{equation} \label{eq:5}
{\mathord{\buildrel{\lower3pt\hbox{$\scriptscriptstyle\rightharpoonup$}}
\over F} _{j + \frac{1}{2}}}({0^ + }) =
\mathord{\buildrel{\lower3pt\hbox{$\scriptscriptstyle\rightharpoonup$}}
\over F} _j^ + ({0^ + }) +
\mathord{\buildrel{\lower3pt\hbox{$\scriptscriptstyle\rightharpoonup$}}
\over F} _{j + 1}^ - ({0^ + }).
\end{equation}

By a first-order Steger-Warming scheme \cite{sw81}, we have
\begin{equation} \label{eq:6}
{\mathord{\buildrel{\lower3pt\hbox{$\scriptscriptstyle\rightharpoonup$}}
\over F} ^ \pm
}(\mathord{\buildrel{\lower3pt\hbox{$\scriptscriptstyle\rightharpoonup$}}
\over W} ) = {A^ \pm
}\mathord{\buildrel{\lower3pt\hbox{$\scriptscriptstyle\rightharpoonup$}} \over W},
\end{equation}
where ${A^ \pm } = L{\Lambda ^ \pm }R$,  $L$ and $R$ are the left and
right eigenvectors respectively, $\Lambda^\pm$ are the diagonal matrices of
positive resp. negative eigenvalues. On the other hand, we also have
\begin{equation} \label{eq:7}
\mathord{\buildrel{\lower3pt\hbox{$\scriptscriptstyle\rightharpoonup$}}
\over F}
(\mathord{\buildrel{\lower3pt\hbox{$\scriptscriptstyle\rightharpoonup$}}
\over W} ) =
A\mathord{\buildrel{\lower3pt\hbox{$\scriptscriptstyle\rightharpoonup$}} \over W},
\end{equation}
where $A = L\Lambda R$ is the Jacobian matrix, $\Lambda =\Lambda^+  +\Lambda^-$
is the diagonal matrix of eigenvalues.

\subsection{Time derivatives}

By virtue of (\ref{eq:6}), (\ref{eq:1}) and (\ref{eq:7}), omitting the subscripts, we obtain
the time derivatives of two parts in (\ref{eq:5}):

\[\frac{\partial }{{\partial t}}{\mathord{\buildrel{\lower3pt\hbox{$\scriptscriptstyle\rightharpoonup$}}
\over F} ^ \pm } = {A^ \pm }\frac{\partial }{{\partial
t}}\mathord{\buildrel{\lower3pt\hbox{$\scriptscriptstyle\rightharpoonup$}}
\over W}  = {A^ \pm }\left( { - A\frac{\partial }{{\partial
x}}\mathord{\buildrel{\lower3pt\hbox{$\scriptscriptstyle\rightharpoonup$}} \over W} } \right).  ???? \]

Generally, we have
\begin{equation} \label{eq:8}
\frac{{{\partial ^k}}}{{\partial
{t^k}}}{\mathord{\buildrel{\lower3pt\hbox{$\scriptscriptstyle\rightharpoonup$}}
\over F} ^ \pm } = \frac{{{\partial ^{k - 1}}}}{{\partial {t^{k -
1}}}}\frac{\partial }{{\partial
t}}{\mathord{\buildrel{\lower3pt\hbox{$\scriptscriptstyle\rightharpoonup$}}
\over F} ^ \pm } = \frac{{{\partial ^{k - 1}}}}{{\partial {t^{k -
1}}}}\left( {{A^ \pm }\left( { - A\frac{\partial }{{\partial
x}}\mathord{\buildrel{\lower3pt\hbox{$\scriptscriptstyle\rightharpoonup$}}
\over W} } \right)} \right) = {A^ \pm }{\left( { - A} \right)^k}\frac{{{\partial ^k}}}{{\partial
{x^k}}}\mathord{\buildrel{\lower3pt\hbox{$\scriptscriptstyle\rightharpoonup$}} \over W}.
\end{equation}
Hence, we can easily get the numerical flux (in high order accuracy) as follows.
\begin{equation} \label{eq:9}
\mathord{\buildrel{\lower3pt\hbox{$\scriptscriptstyle\rightharpoonup$}}
\over F} _{j + \frac{1}{2}} (\tau ) =
\mathord{\buildrel{\lower3pt\hbox{$\scriptscriptstyle\rightharpoonup$}}
\over F} _{j + \frac{1}{2}} (0^ +  ) + A^ \pm  \sum\limits_{k = 1}^N
{\left( { - A} \right)^k \frac{{\partial ^k
\mathord{\buildrel{\lower3pt\hbox{$\scriptscriptstyle\rightharpoonup$}}
\over W} }}{{\partial x^k }}\frac{{\tau ^k }}{{k!}}},
 \end{equation}
 where the spatial derivatives $\partial^k\mathord{\buildrel{\lower3pt\hbox{$\scriptscriptstyle\rightharpoonup$}}\over W}/\partial x^k$
 ($k=1,\cdots, N$) are given in \ref{eq:21}-\ref{eq:23}.

\textbf{Remark 3.} From the above construction procedure, we can see
the new scheme HFVS has two advantages. Firstly , it is obviously
simpler than the original ADER and HGKS. The main reason is: the
high order time derivatives can be expressed by a explicit and
compact formulation of high order space derivatives, which makes the
computational procedures clear and easy to implement. In addition,
HFVS evidently used less computational costs than ADER and HGKS.
This is also attributed to the simpler procedure to compute the high
derivative terms.

\textbf{Remark 4.} Although the Steger-Warming scheme is used as the
leading term in the construction of our scheme,  we should point out
that the leading term in (\ref{eq:9}) can be arbitrary Riemann
solver, such as the HLLC and BGK solvers. This means the latter
terms can be used as a building block to extend the current FVS
scheme to high order accuracy in both space and time.

\section{Numerical tests}
In this section, we present the numerical experiment results for
linear and nonlinear hyperbolic conservative laws in 1D and 2D. We
call our new scheme the HFVS (arbitrary High order Flux Vector
Splitting) scheme. The abbreviations HFVS2, HFVS3 and HFVS5 mean the
second-, third- and fifth-order HFVS schemes, respectively. owing to
the simplicity, all numerical schemes only use conservative
variables ENO/WENO reconstruction techniques. And we will find they
work well enough in most of the cases below.

\subsection{Linear hyperbolic problems}
In this subsection, we want to examine the accuracy of the new schemes by solving some linear hyperbolic problems.

 \textbf{Example 1. Accuracy test}.

 We consider the following linear convection equation
\[\frac{{\partial w}}{{\partial t}} + \frac{{\partial w}}{{\partial x}} = 0\]
with initial data
\[{w_0}(x) = \sin (2\pi x) .\]
We take the unit interval $(0,1)$ as the computational region and prescribe suitable boundary conditions (inflow boundary
condition).

Table I shows the errors and accuracy for different schemes using a
same $CFL$ number $0.95$ at a computational time $t=1.0$. We can
observe that all HFVS schemes can reach their designed order
accuracy.


\begin{table}[htbp]
 \begin{center}   \caption{Errors and accuracy for WENO and HFVS}
 \label{table1}
   \begin{tabular}{|c|c|c|c|c|c|c|c|}
     \hline
     METHOD & $N$& $L_1 ERROR$ &ORDER& $L_2 ERROR$ &ORDER& $L_{\infty} ERROR$ &ORDER
     \\
     \hline

       HFVS2 &20    &0.28    &       &0.28    &         &0.31    &  \\
             &40    &0.57E-01&  2.34 &0.65E-01&  2.15   &0.92E-01&  1.76\\
             &80    &0.12E-01&  2.21 &0.13E-01&  2.25   &0.19E-01&  2.24\\
            &160    &0.29E-02&  2.07 &0.32E-02&  2.08   &0.46E-02&  2.08\\
            &320    &0.73E-03&  2.01 &0.81E-03&  2.01   &0.11E-02&  2.01\\
            &640    &0.18E-03&  2.01 &0.20E-03&  2.01   &0.28E-03&  2.01\\
     \hline

       HFVS3 &20    &2.32E-01&       &2.54E-01&              &3.35E-01&  \\
             &40    &3.92E-02&  2.56 &4.38E-02&  2.53   &6.14E-02&  2.44\\
             &80    &5.26E-03&  2.89 &5.86E-03&  2.90   &8.27E-03&  2.89\\
            &160    &6.62E-04&  2.99 &7.36E-04&  2.99   &1.04E-03&  2.99\\
            &320    &8.13E-05&  3.02 &9.03E-05&  3.02   &1.28E-04&  3.02\\
            &640    &1.02E-05&  2.99 &1.13E-05&  2.99   &1.60E-05&  3.00\\
       \hline
       HFVS5 &20    &1.02E-01&                &1.17E-01&              &1.59E-01&  \\
             &40    &4.75E-03&  4.42 &5.35E-03&  4.45 &7.51E-03&  4.40\\
             &80    &1.62E-04&  4.87 &1.81E-04&  4.88 &2.55E-04&  4.88\\
            &160    &5.14E-06&  4.97 &5.71E-06&  4.98 &8.08E-06&  4.97\\
            &320    &1.58E-07&  5.02 &1.75E-07&  5.02 &2.48E-07&  5.02\\
            &640    &4.93E-09&  5.00 &5.48E-09&  4.99 &7.75E-09&  5.00\\
       \hline
   \end{tabular}
 \end{center}
 \end{table}

\subsection{Nonlinear hyperbolic problem}

In this section, we continue to present the numerical results of 1D
and 2D Euler equations of compressible flows.

\subsubsection{1D Euler equations}
The Governing equations of 1D Euler equations can see \ref{1deuler}.
In this section, $\gamma=1.4$ ,and all numerical schemes use an
optimal $CFL$ number $0.95$.

\textbf{Example 2. Shu-Osher shock turbulence interaction
problem\cite{shu-osher}}.

 the computational domain is $[-1,1]$, and the
initial condition is
\[\left( {\rho ,u,P} \right) = \left\{ {\begin{array}{*{20}{c}}
   {(3.857143,2.629369,10.333333),x <  - 0.8}  \\
   {(1 + 0.2\sin 5\pi x,0.0,1.0) \quad \quad \quad \;,x >  - 0.8}  \\
\end{array}} \right.\]

\begin{figure}[!th]
\renewcommand{\baselinestretch}{0.01}
\setlength{\unitlength}{1mm} \setlength{\fboxsep}{0pt}
\begin{center}
\scalebox{0.35}{\includegraphics*[80,60][650,550]{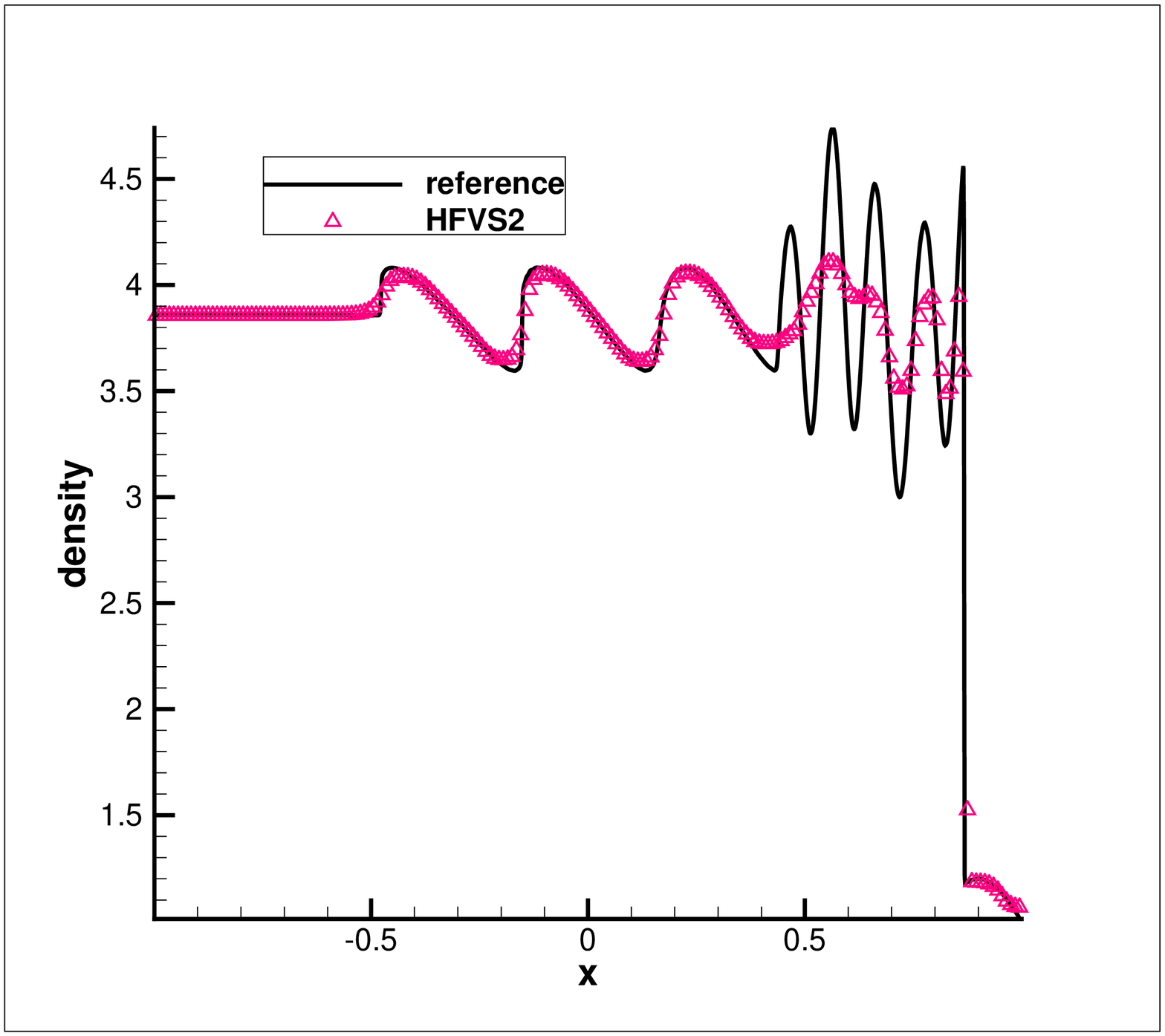}}\quad
\scalebox{0.35}{\includegraphics*[80,60][650,550]{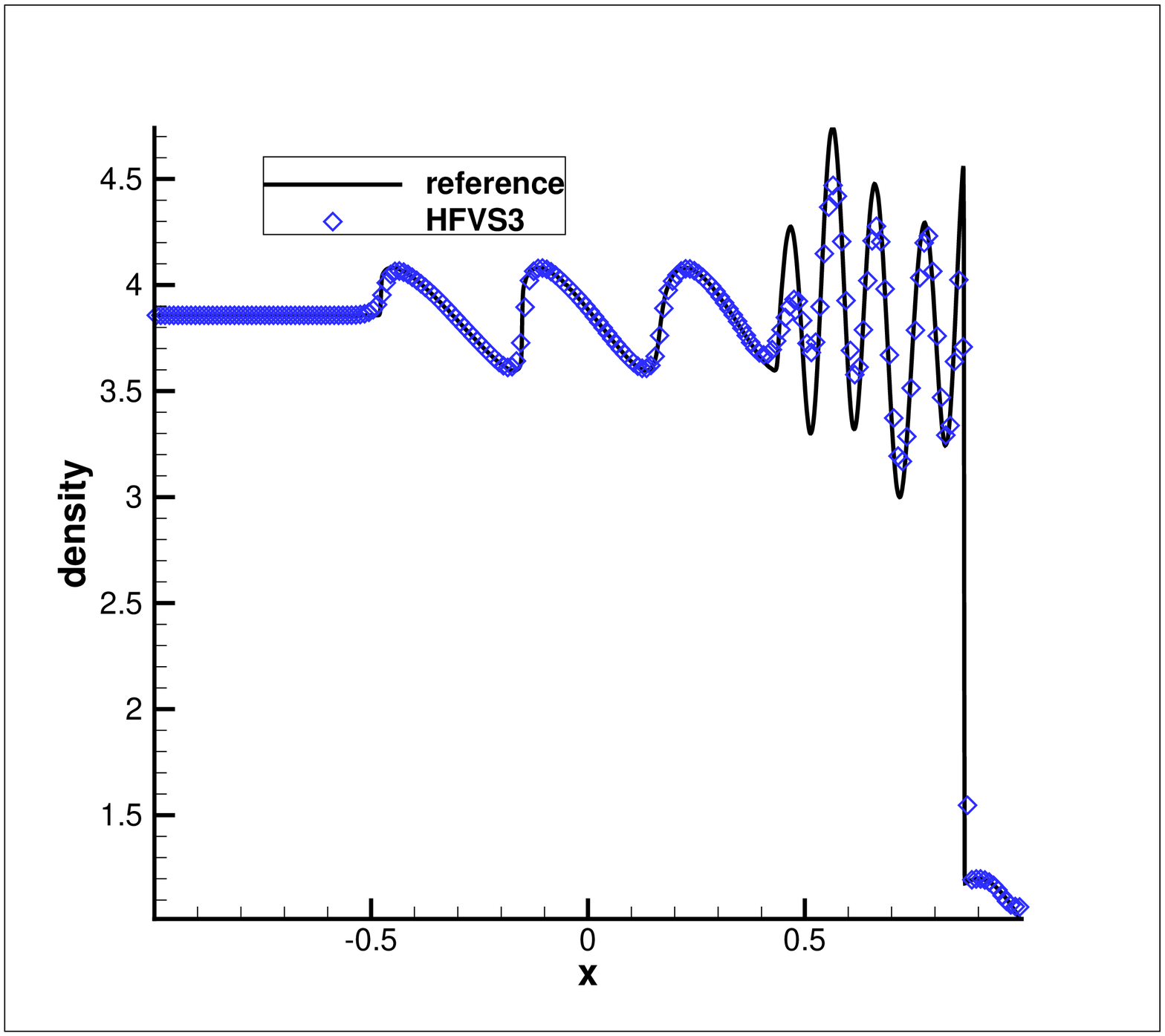}}\quad
\scalebox{0.35}{\includegraphics*[80,60][650,550]{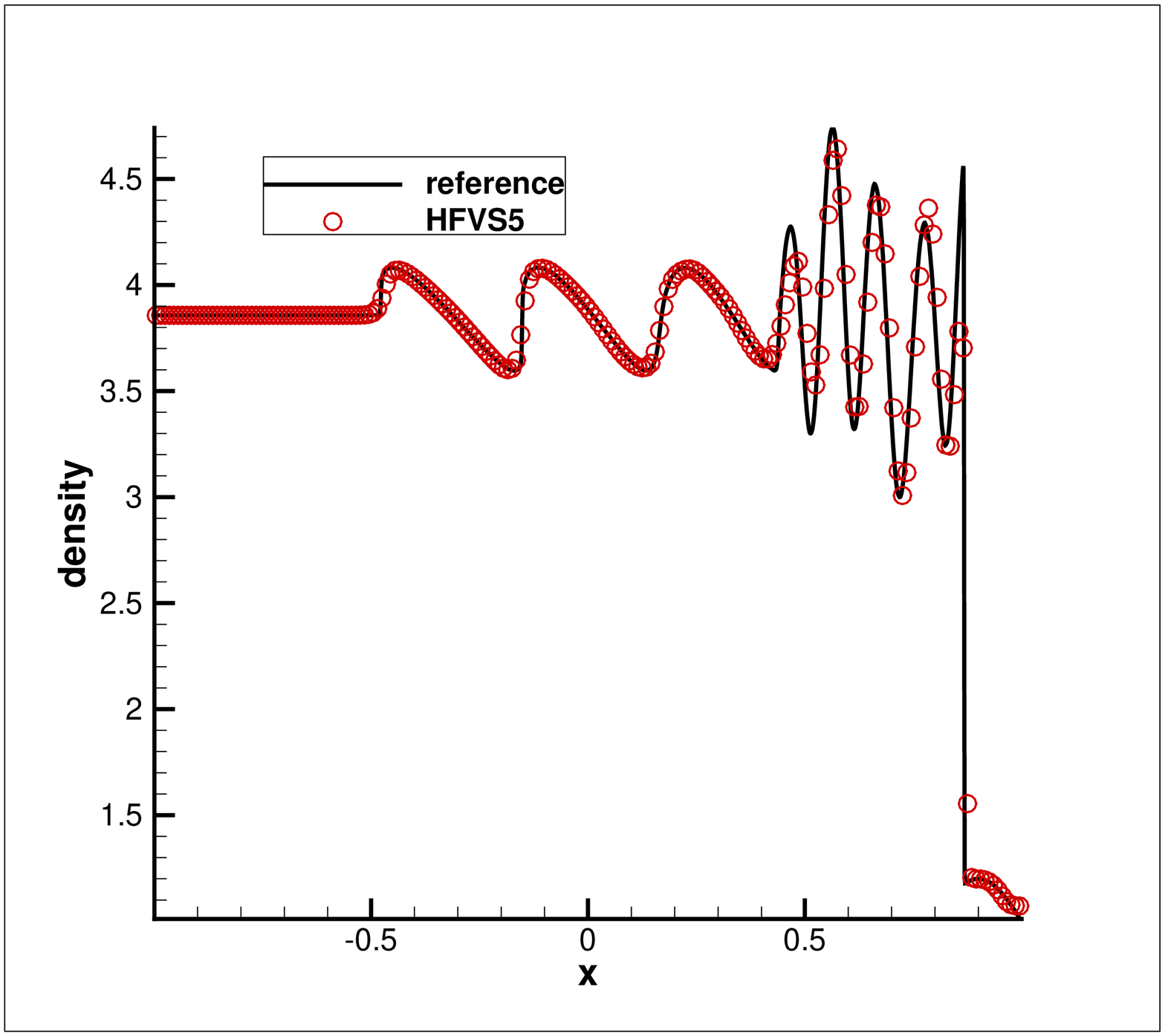}}
\end{center}
\caption{ Example 2. density profiles for different
schemes.\label{Exam2.1}}
\end{figure}

\begin{figure}[!th]
\renewcommand{\baselinestretch}{0.01}
\setlength{\unitlength}{1mm} \setlength{\fboxsep}{0pt}
\begin{center}
\scalebox{0.35}{\includegraphics*[80,60][650,550]{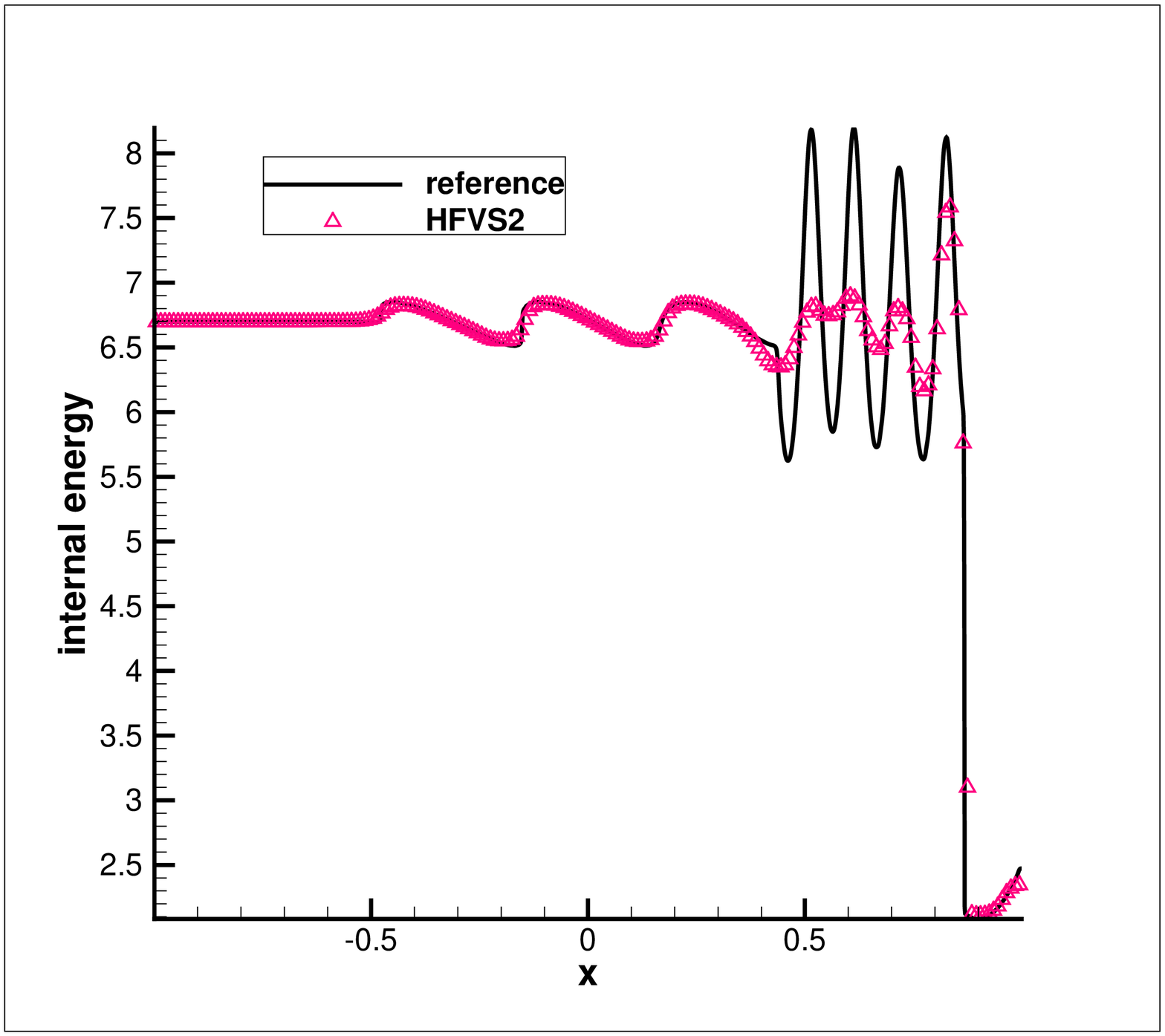}}\quad
\scalebox{0.35}{\includegraphics*[80,60][650,550]{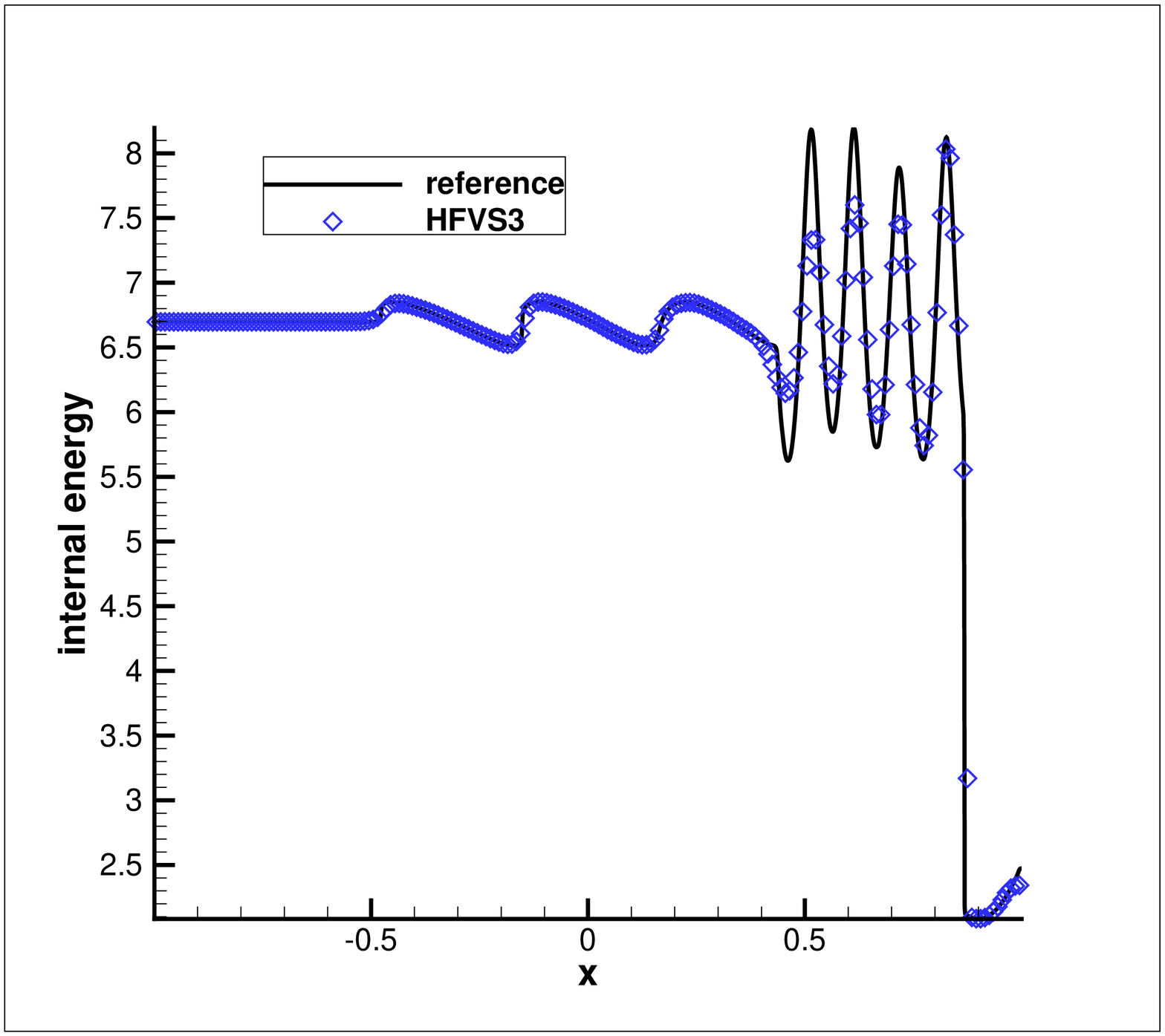}}\quad
\scalebox{0.35}{\includegraphics*[80,60][650,550]{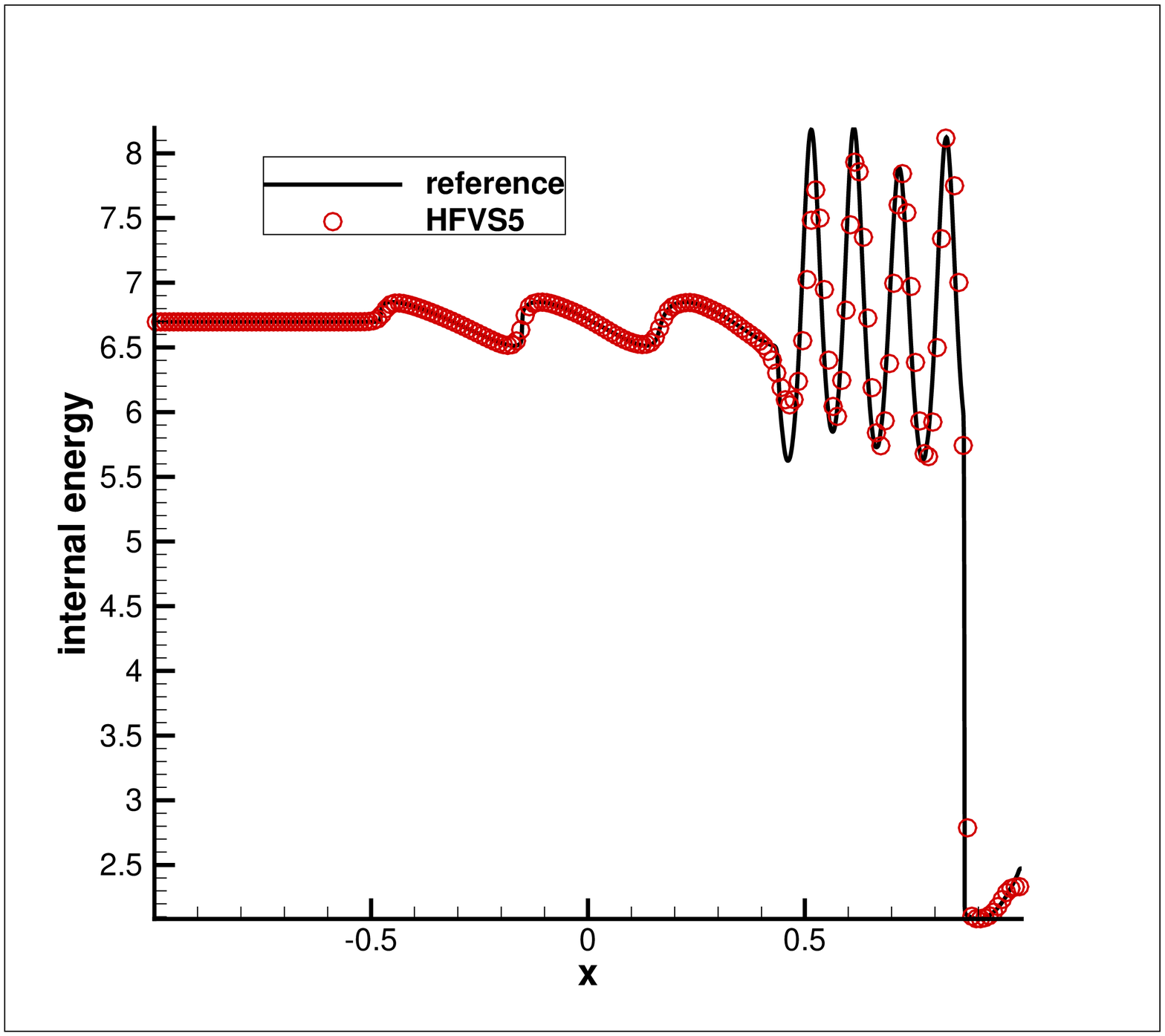}}
\end{center}
\caption{ Example 2. internal energy profiles for different
schemes.\label{Exam2.2}}
\end{figure}

Figs.\ref{Exam2.1} and Figs.\ref{Exam2.2} shows the density ,
pressure and internal profiles computed by different schemes
$t=0.47$ with mesh of 200 cells. And a reference solution is
computed by the third order WENO scheme on a fine mesh with 2000
cells. We can observe that all schemes can resolve the complex
solution of the governing equation. In addition, HFVS2 introduces
over diffusion, and HFVS3 improves the resolution, while HFVS5 gives
a sharp profile. This also demonstrates the advantages of high order
scheme.

\textbf{Example 3 Woodward-Colella blast wave
problem\cite{Woodward84}}. The computational domain is $[0,1]$ with
reflected boundary conditions on both sides, and the initial
condition is

\[\left( {\rho ,U,P} \right) = \left\{ {\begin{array}{*{20}{c}}
   {(1,0,1000),0 \le x \le 0.1}  \\
   {(1,0,0.01),0.1 \le x \le 0.9}  \\
   {(1,0,100),0.9 \le x \le 1.0}  \\
\end{array}} \right.\]

This is problem with high ratio of pressure and the solution
includes strong shock waves, contact discontinuities, rarefaction
waves and their interactions. So it is a very challenging problem
for high order schemes.

\begin{figure}[!th]
\renewcommand{\baselinestretch}{0.01}
\setlength{\unitlength}{1mm} \setlength{\fboxsep}{0pt}
\begin{center}
\scalebox{0.35}{\includegraphics*[80,60][650,550]{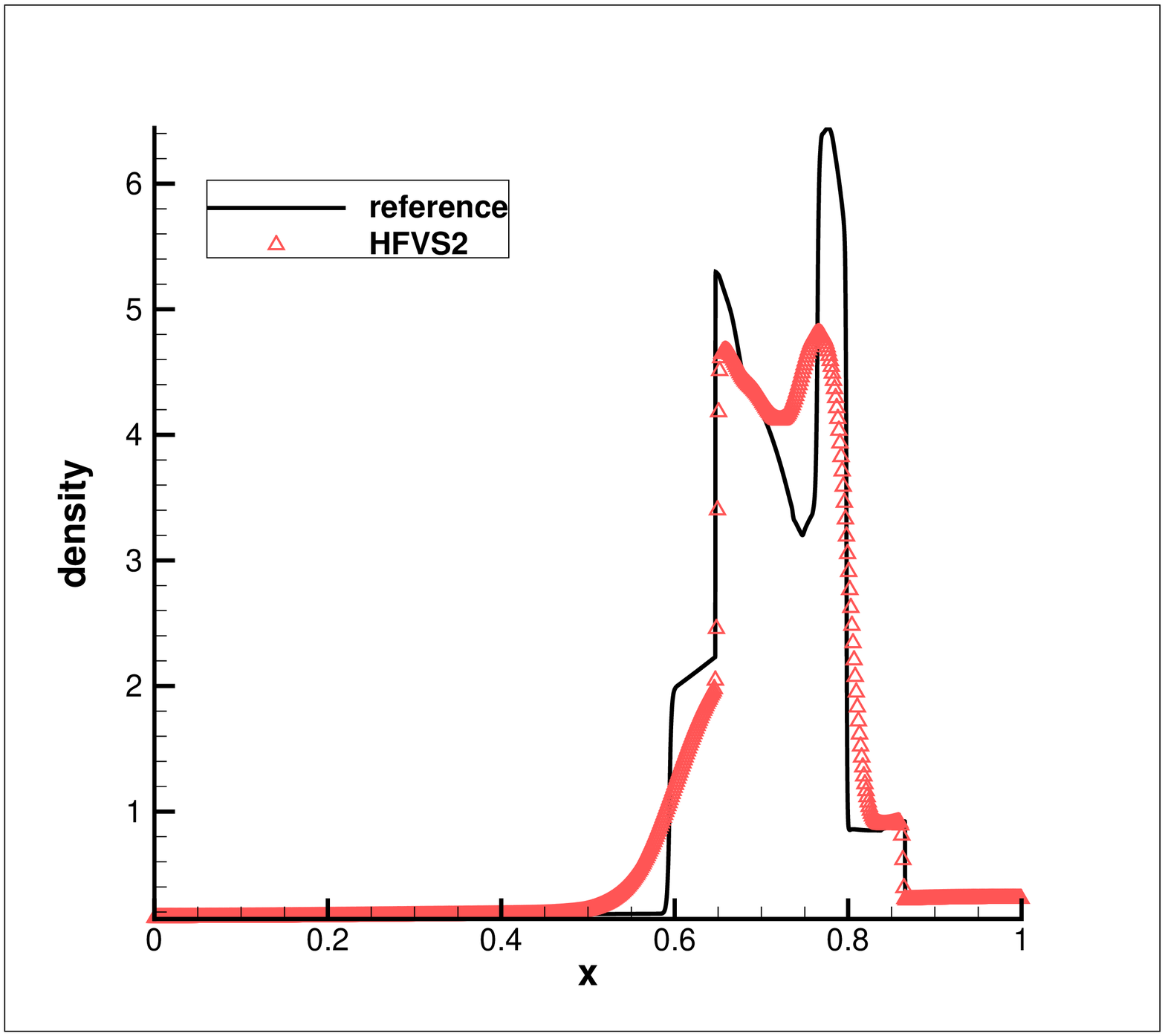}}\quad
\scalebox{0.35}{\includegraphics*[80,60][650,550]{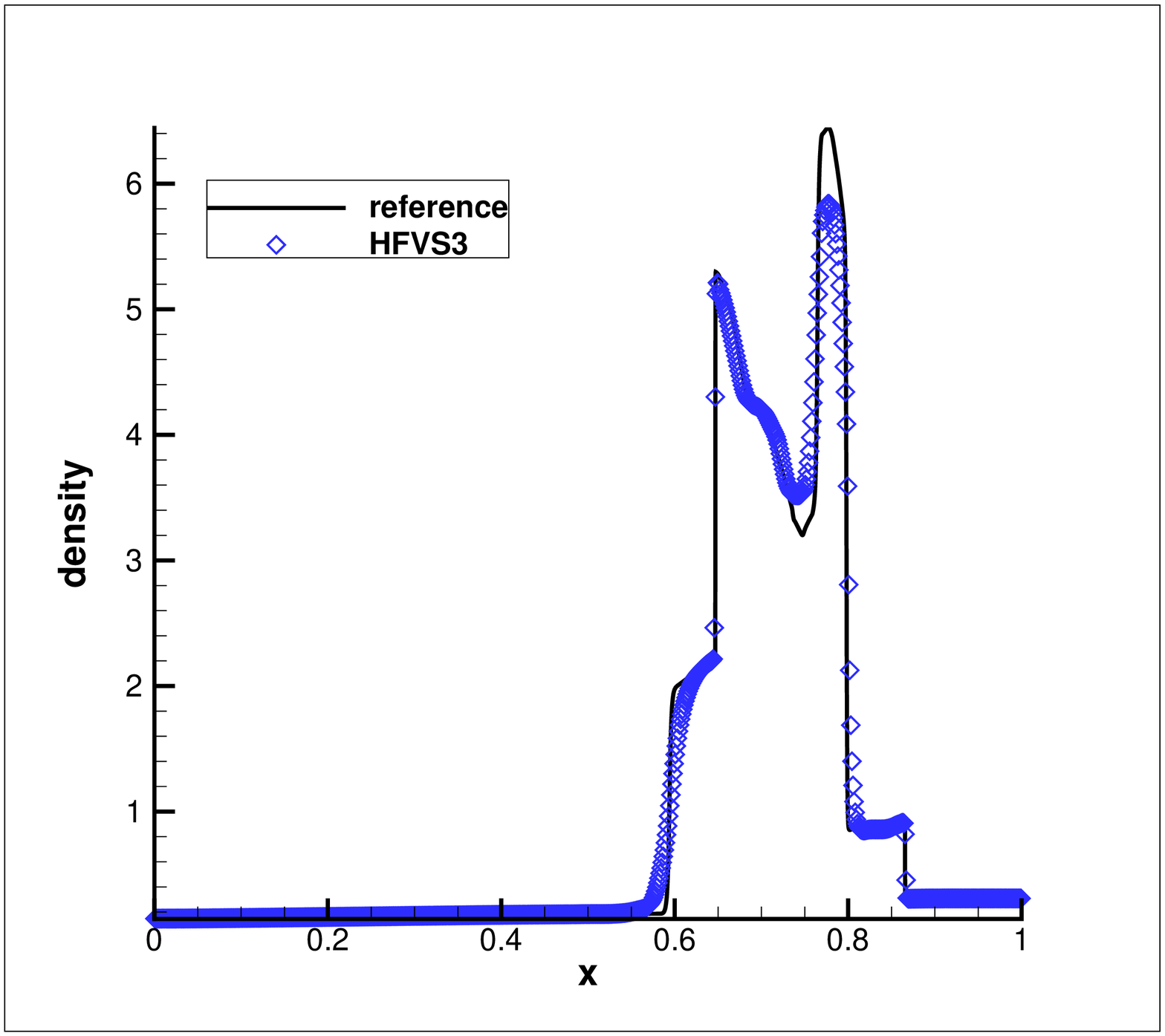}}\quad
\scalebox{0.35}{\includegraphics*[80,60][650,550]{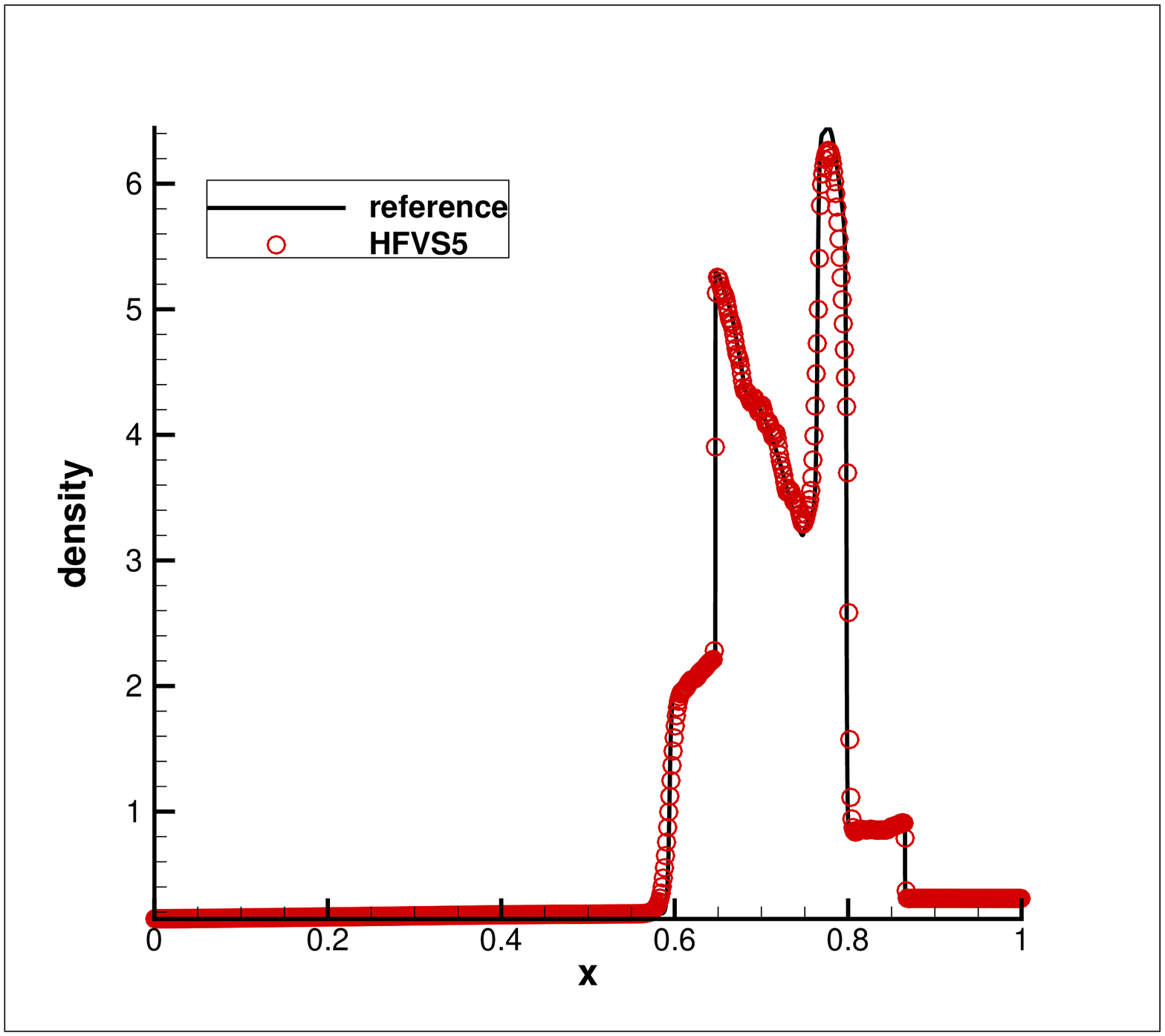}}
\end{center}
\caption{ Example 3. density profiles for different
schemes.\label{Exam3.1}}
\end{figure}

\begin{figure}[!th]
\renewcommand{\baselinestretch}{0.01}
\setlength{\unitlength}{1mm} \setlength{\fboxsep}{0pt}
\begin{center}
\scalebox{0.35}{\includegraphics*[80,60][650,550]{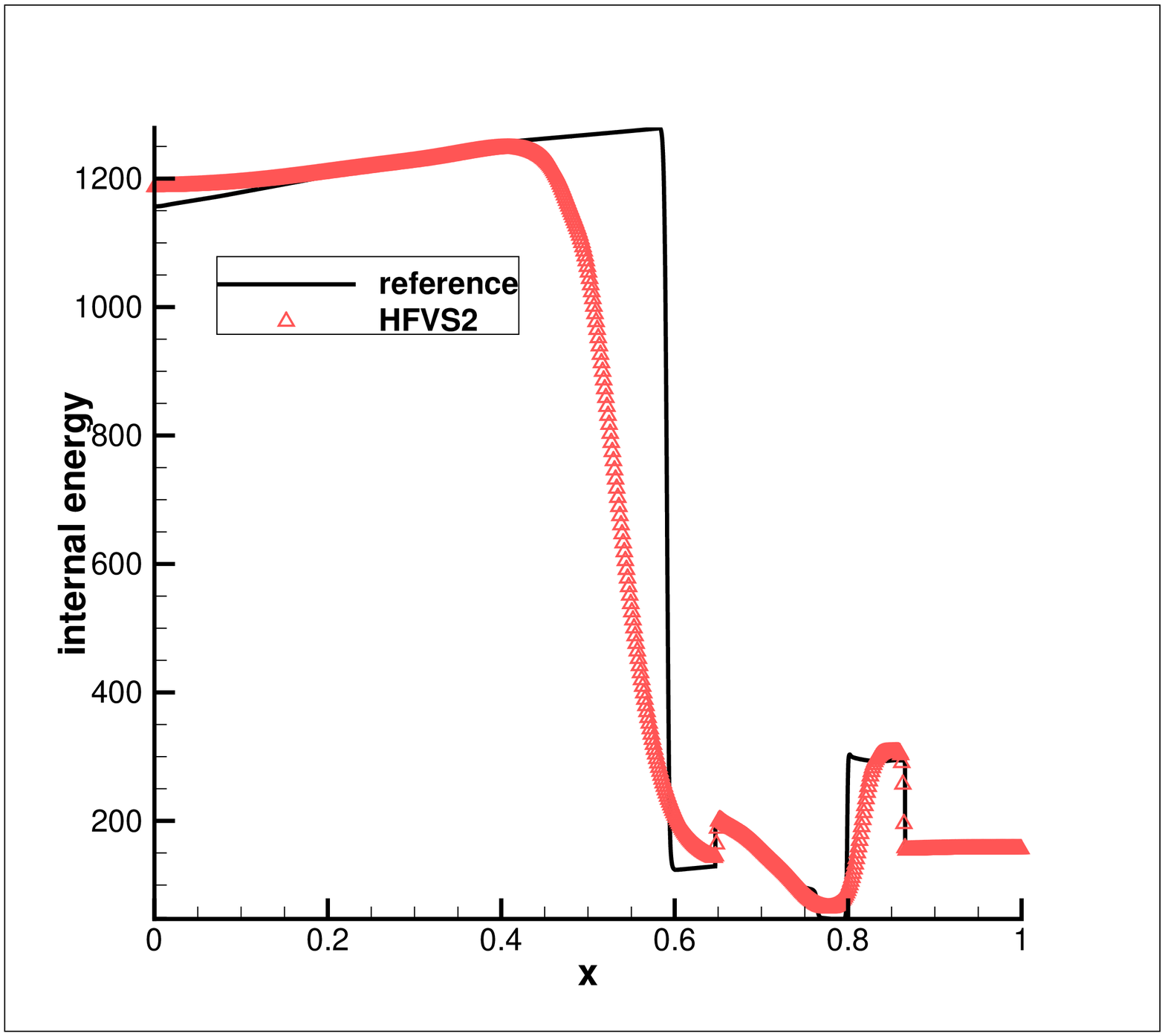}}\quad
\scalebox{0.35}{\includegraphics*[80,60][650,550]{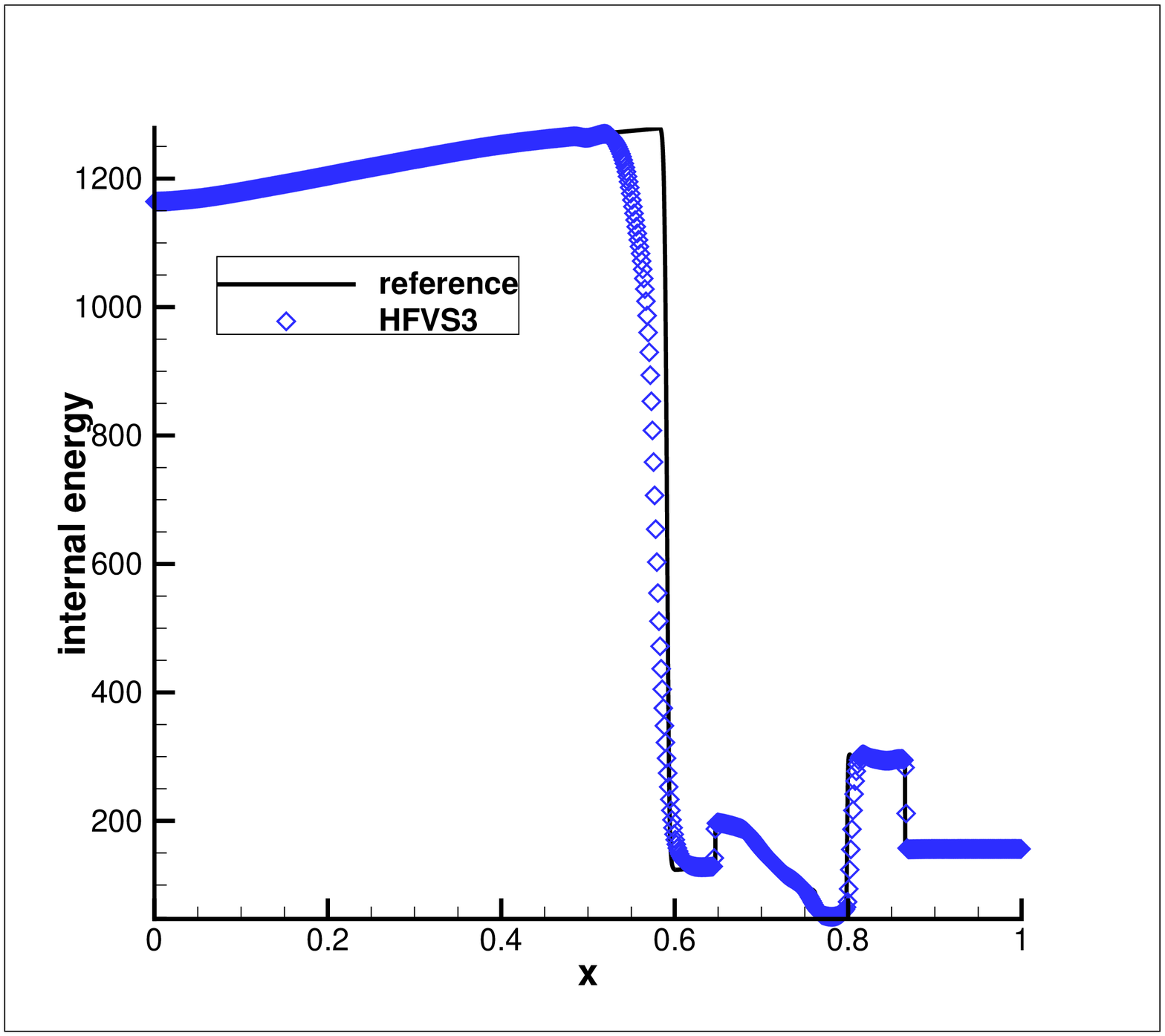}}\quad
\scalebox{0.35}{\includegraphics*[80,60][650,550]{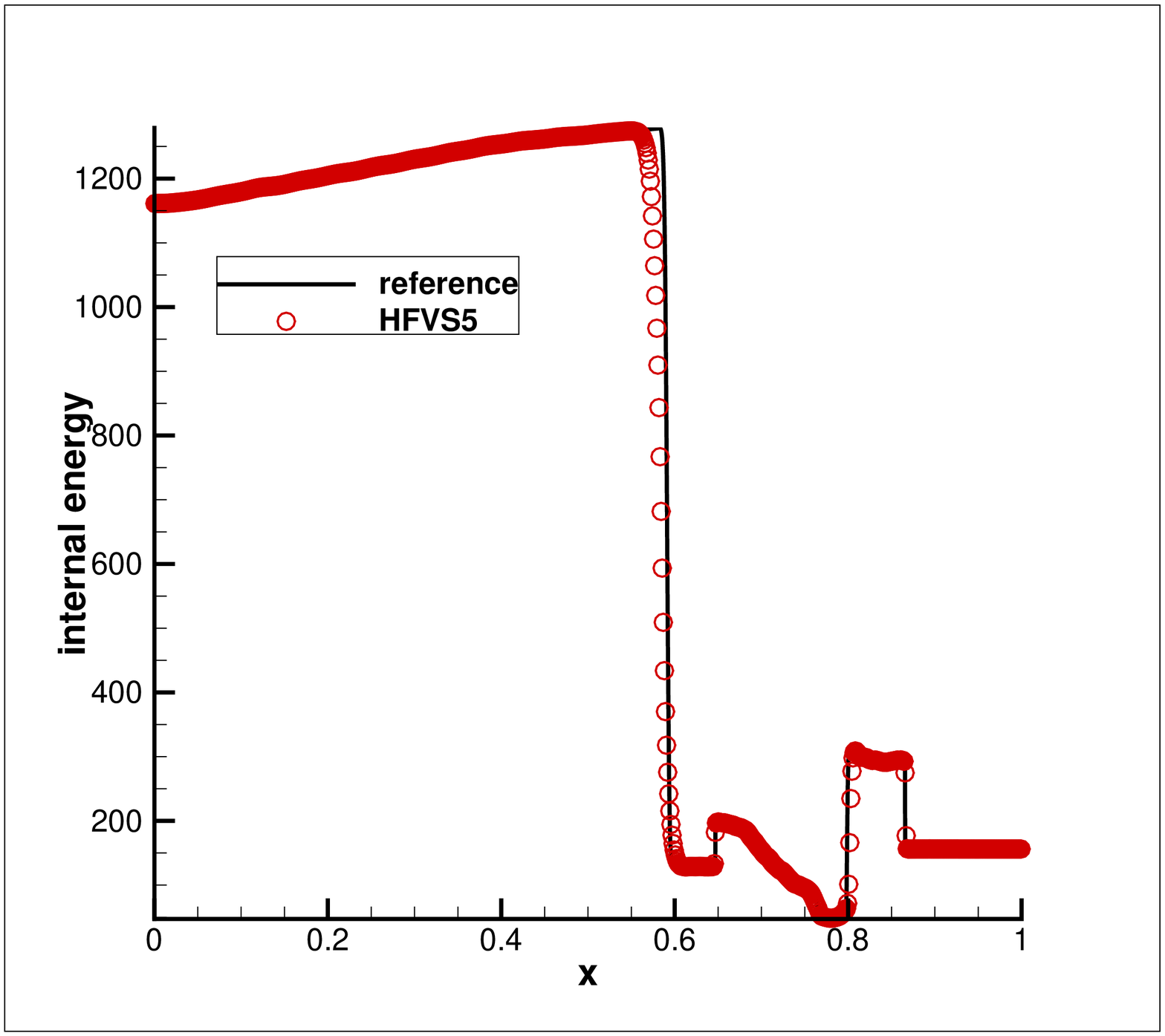}}
\end{center}
\caption{ Example 3. internal energy profiles for different
schemes.\label{Exam3.2}}
\end{figure}

 Fig.\ref{Exam3.1} and Fig.\ref{Exam3.2}  shows the numerical
results at computational time $t=0.038$ with mesh of 800 cells. The
reference solution is computed by third order WENO on 10000 cells.
All HFVS schemes can resolve all the waves of the flow. As the
previous test problem, HFVS2 gives poor solution, HFVS3 greatly
improve the resolution, and HFVS5 can give a sharp picture. The
numerical results can demonstrate the robustness and accuracy of
HFVS.

\subsubsection{2D nonlinear problems} In this section, we present the
numerical results of 2D Euler equations:
\[
\frac{{\partial
\mathord{\buildrel{\lower3pt\hbox{$\scriptscriptstyle\rightharpoonup$}}
\over W} }}{{\partial t}} + \frac{{\partial
\mathord{\buildrel{\lower3pt\hbox{$\scriptscriptstyle\rightharpoonup$}}
\over F}
(\mathord{\buildrel{\lower3pt\hbox{$\scriptscriptstyle\rightharpoonup$}}
\over W} )}}{{\partial x}} + \frac{{\partial
\mathord{\buildrel{\lower3pt\hbox{$\scriptscriptstyle\rightharpoonup$}}
\over G}
(\mathord{\buildrel{\lower3pt\hbox{$\scriptscriptstyle\rightharpoonup$}}
\over W} )}}{{\partial x}} = 0
\]
where
\[
\mathord{\buildrel{\lower3pt\hbox{$\scriptscriptstyle\rightharpoonup$}}
\over W}  = \left( {\rho ,\rho U,\rho V,\rho E} \right)^T
\]

\[
\mathord{\buildrel{\lower3pt\hbox{$\scriptscriptstyle\rightharpoonup$}}
\over F}
(\mathord{\buildrel{\lower3pt\hbox{$\scriptscriptstyle\rightharpoonup$}}
\over W} ) = \left( {\rho U,\rho U^2  + P,\rho UV,\rho EU + PU}
\right)^T
\]
\[
\mathord{\buildrel{\lower3pt\hbox{$\scriptscriptstyle\rightharpoonup$}}
\over G}
(\mathord{\buildrel{\lower3pt\hbox{$\scriptscriptstyle\rightharpoonup$}}
\over W} ) = \left( {\rho V,\rho UV,\rho V^2  + P,\rho EV + PV}
\right)^T
\]
and the equations of state is
\[
P = \left( {\gamma  - 1} \right)\rho \left( {E - \frac{1}{2}U^2  -
\frac{1}{2}V^2 } \right)
\]
In the following case , $\gamma=1.4$.

In this section, all numerical schemes use an optimal $CFL$ number
of $0.45$.

\textbf{ Example 4. 2D Riemann problem}. We want to show two cases
which are simple extensions of 1D Riemann problems.

The first case sets computational domain as $ [0,1] \times [0,1] $.
And the boundary conditions are all outflow condition. The initial
condition is given as

\[
\left( {\rho ,U,V,P} \right) = \left\{ {\begin{array}{*{20}c}
   {(1.0,0.0,0.0,1.0)\;,\;x < 0.5\;and\;y < 0.5}  \\
   {(0.1,0.0,0.0,0.1)\;,\;x < 0.5\;and\;y > 0.5}  \\
   {(1.0,0.0,0.0,1.0),\;x > 0.5\;and\;y > 0.5}  \\
   {(0.1,0.0,0.0,0.1),\;x > 0.5\;and\;y < 0.5}  \\
\end{array}} \right.
\]

We continue to compare the numerical results of HFVS with different
order accuracy. The numerical results using the computational cells
of $100 \times 100$ at $t=0.15$ are presented. Fig.\ref{Exam4.1}
shows the density contours between the range of 0.1 and 1.0. We can
draw the same conclusions as 1D case, i.e., all schemes can work
well and HFVS can improve the accuracy
 and resolution if we increase the order of accuracy.

\begin{figure}[!th]
\renewcommand{\baselinestretch}{0.01}
\setlength{\unitlength}{1mm} \setlength{\fboxsep}{0pt}
\begin{center}
\scalebox{0.35}{\includegraphics*[80,60][650,550]{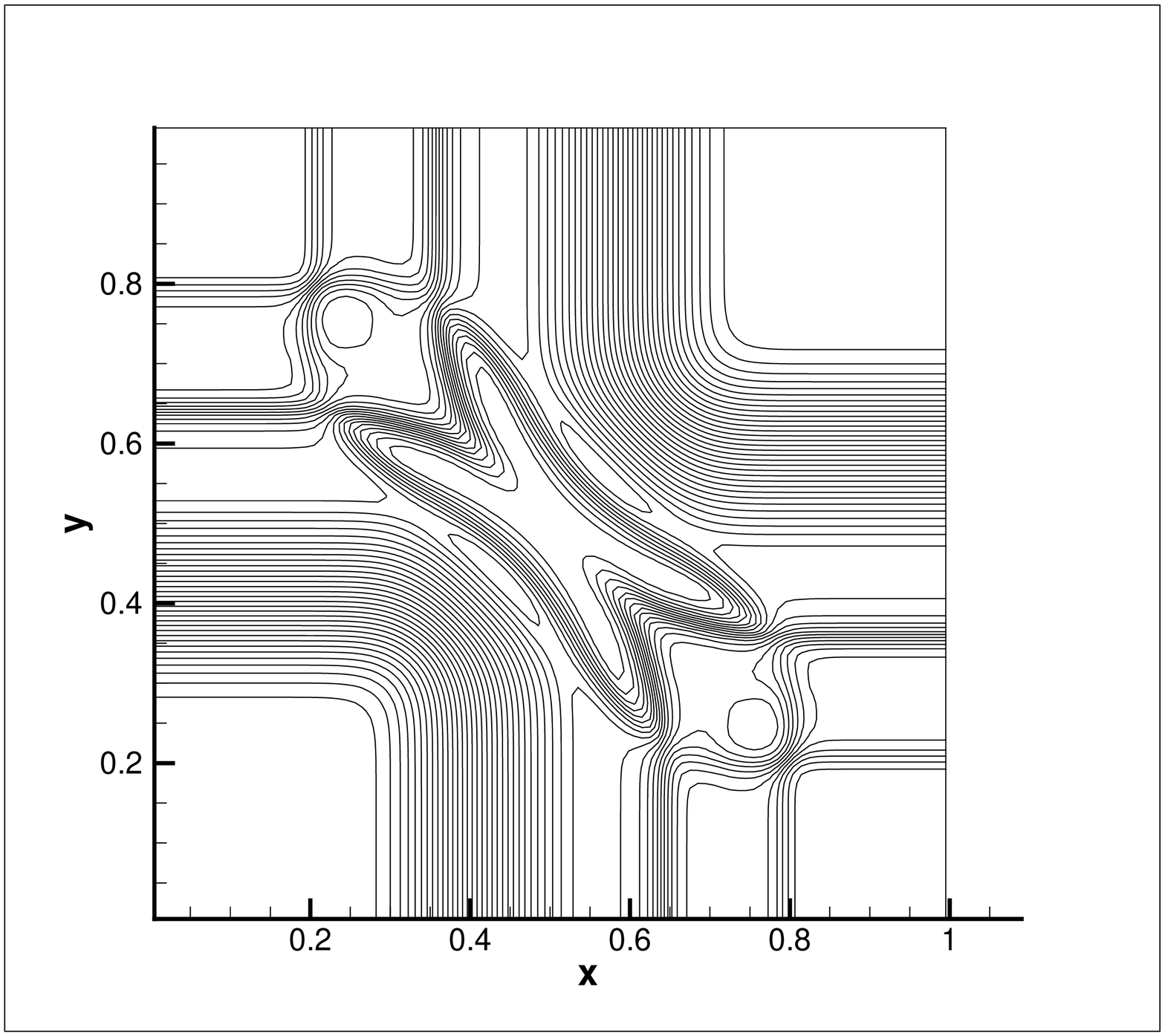}}\quad
\scalebox{0.35}{\includegraphics*[80,60][650,550]{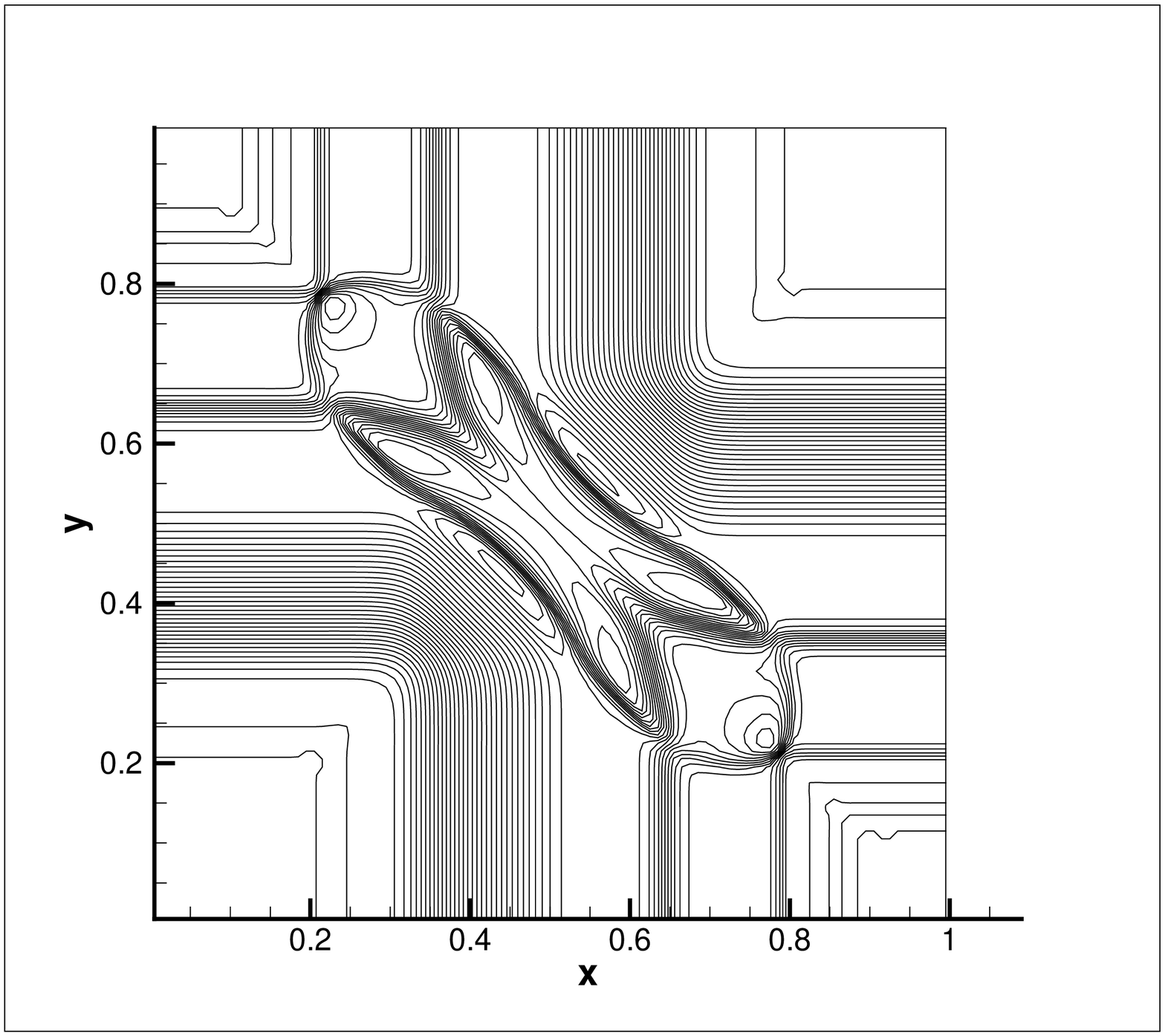}}\quad
\scalebox{0.35}{\includegraphics*[80,60][650,550]{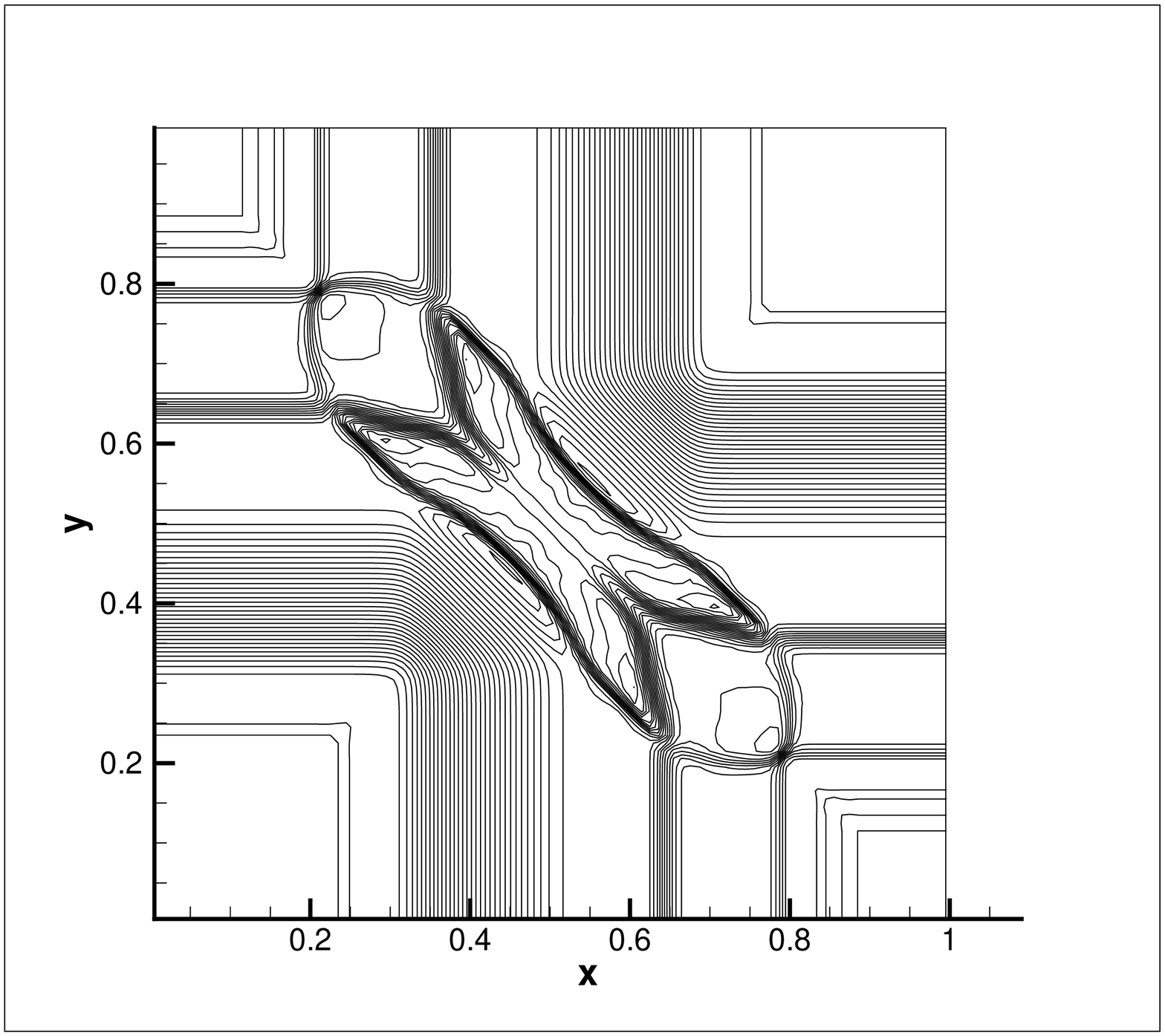}}
\end{center}
\caption{ Example 4. density contours for different schemes,from
left to right: HFVS2, HFVS3, HFVS5. 50 contours were fit between a
range of 0.1 to 1.0.\label{Exam4.1}}
\end{figure}

The second case is a long time test. The computational domain is $
[0,1] \times [0,1] $. And the boundary conditions are all reflective
conditions. The initial condition is given as

\[
\left( {\rho ,U,V,P} \right) = \left\{ {\begin{array}{*{20}c}
   {(1,0,0,1),\sqrt {\left( {x - 0.5} \right)^2  + \left( {y - 0.5} \right)^2 }  \le 0.3}  \\
   {\left( {0.125,0,0,0.1} \right),else}  \\
\end{array}} \right.
\]

The numerical results using the computational cells of $100 \times
100$ at $t=1.0$ are presented. Fig.\ref{Exam4.1} shows the density
contours between the range of 0.1 and 0.6. We can draw the same
conclusions as previous one, i.e., second order scheme gives a
over-diffusion picture, third order scheme improve the accuracy, and
fifth order scheme can give a sharp profile. It also demonstrates
the efficiency of high order schemes.

\begin{figure}[!th]
\renewcommand{\baselinestretch}{0.01}
\setlength{\unitlength}{1mm} \setlength{\fboxsep}{0pt}
\begin{center}
\scalebox{0.35}{\includegraphics*[80,60][650,550]{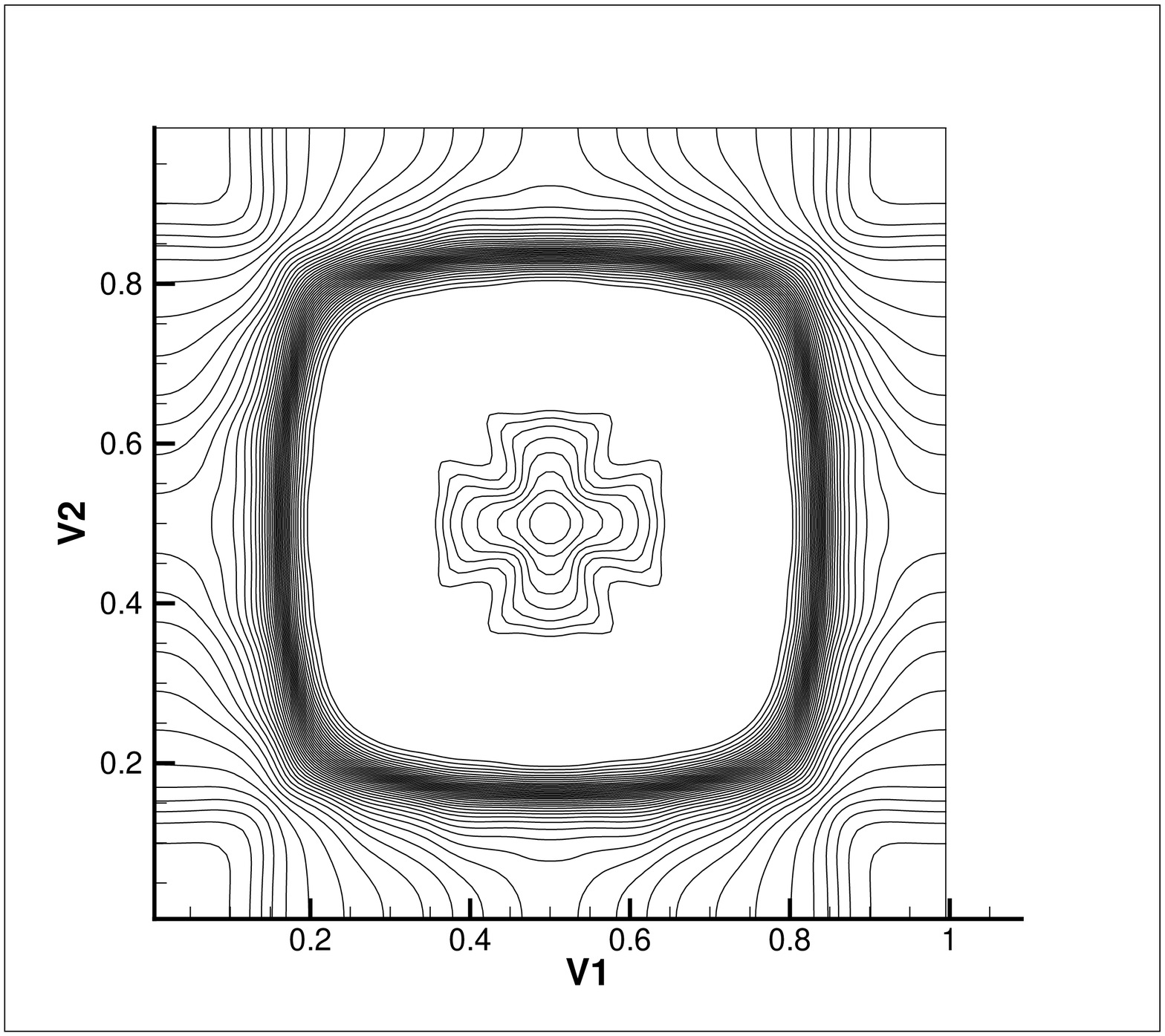}}\quad
\scalebox{0.35}{\includegraphics*[80,60][650,550]{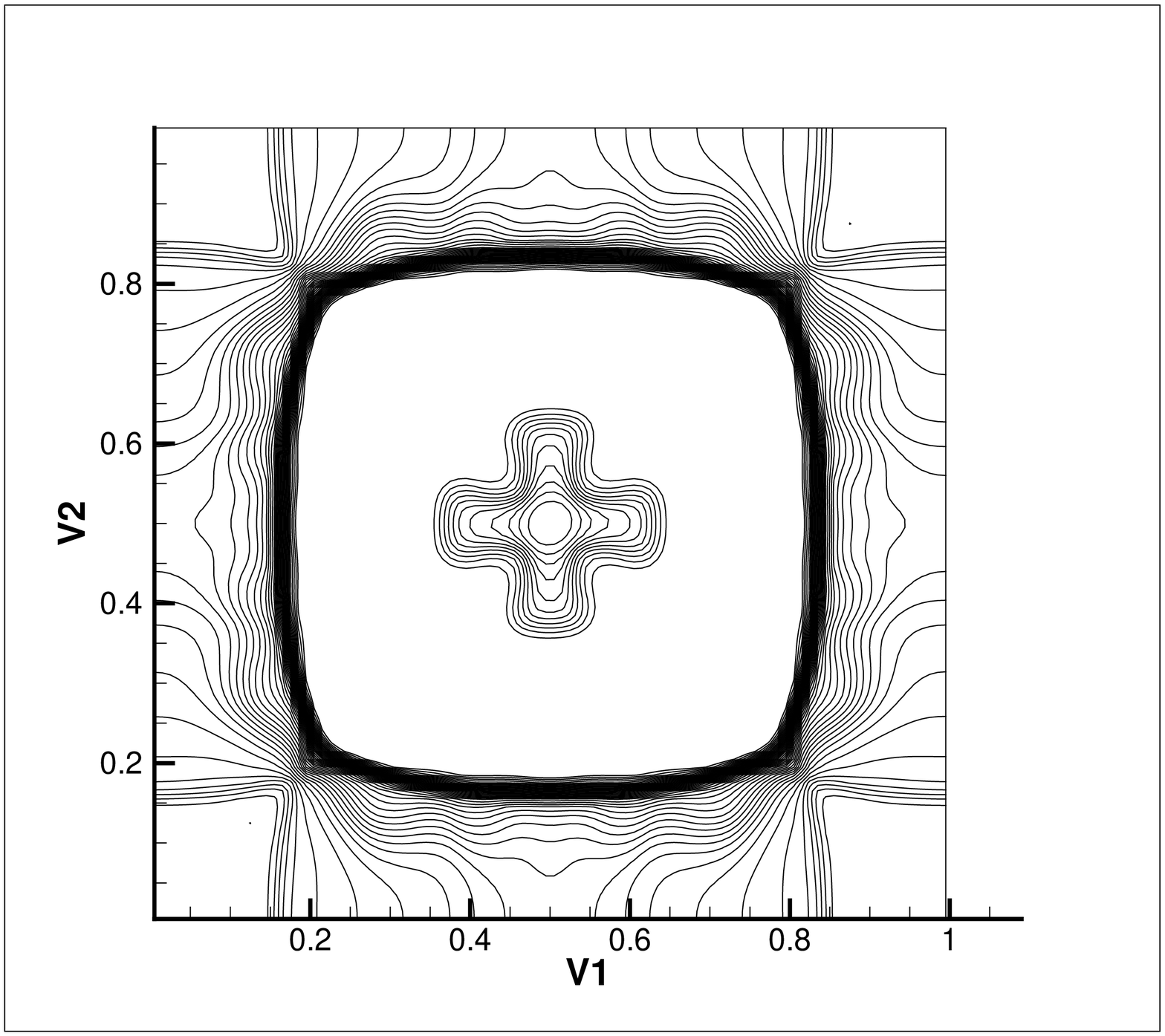}}\quad
\scalebox{0.35}{\includegraphics*[80,60][650,550]{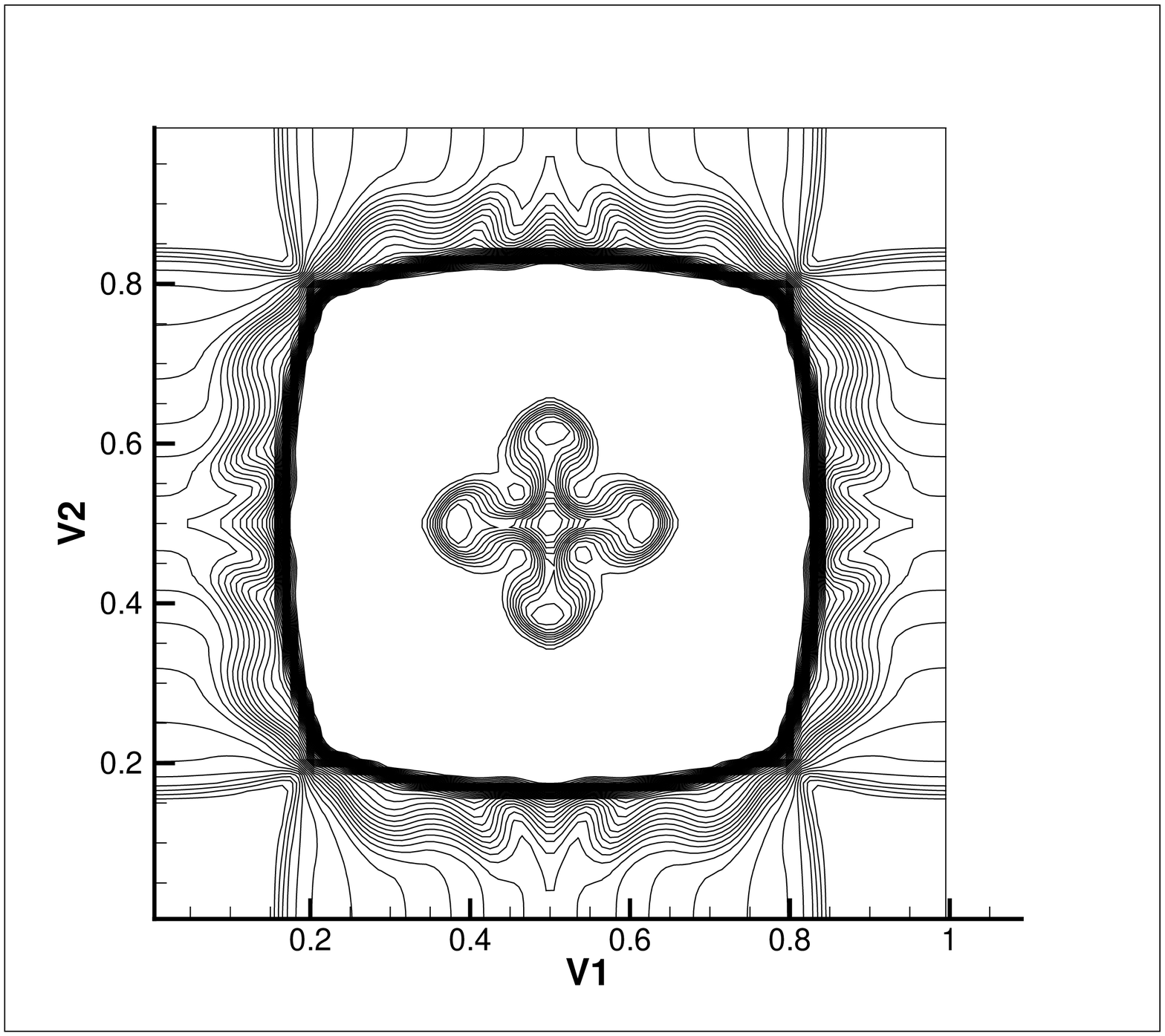}}
\end{center}
\caption{ Example 4. density contours for different schemes,from
left to right: HFVS2, HFVS3, HFVS5. 50 contours were fit between a
range of 0.1 to 1.0.\label{Exam4.1}}
\end{figure}

\subsection{Further discussions}

\subsubsection{The effective of different leading terms}
\begin{figure}[!th]
\renewcommand{\baselinestretch}{0.01}
\setlength{\unitlength}{1mm} \setlength{\fboxsep}{0pt}
\begin{center}

\scalebox{0.35}{\includegraphics*[80,60][650,550]{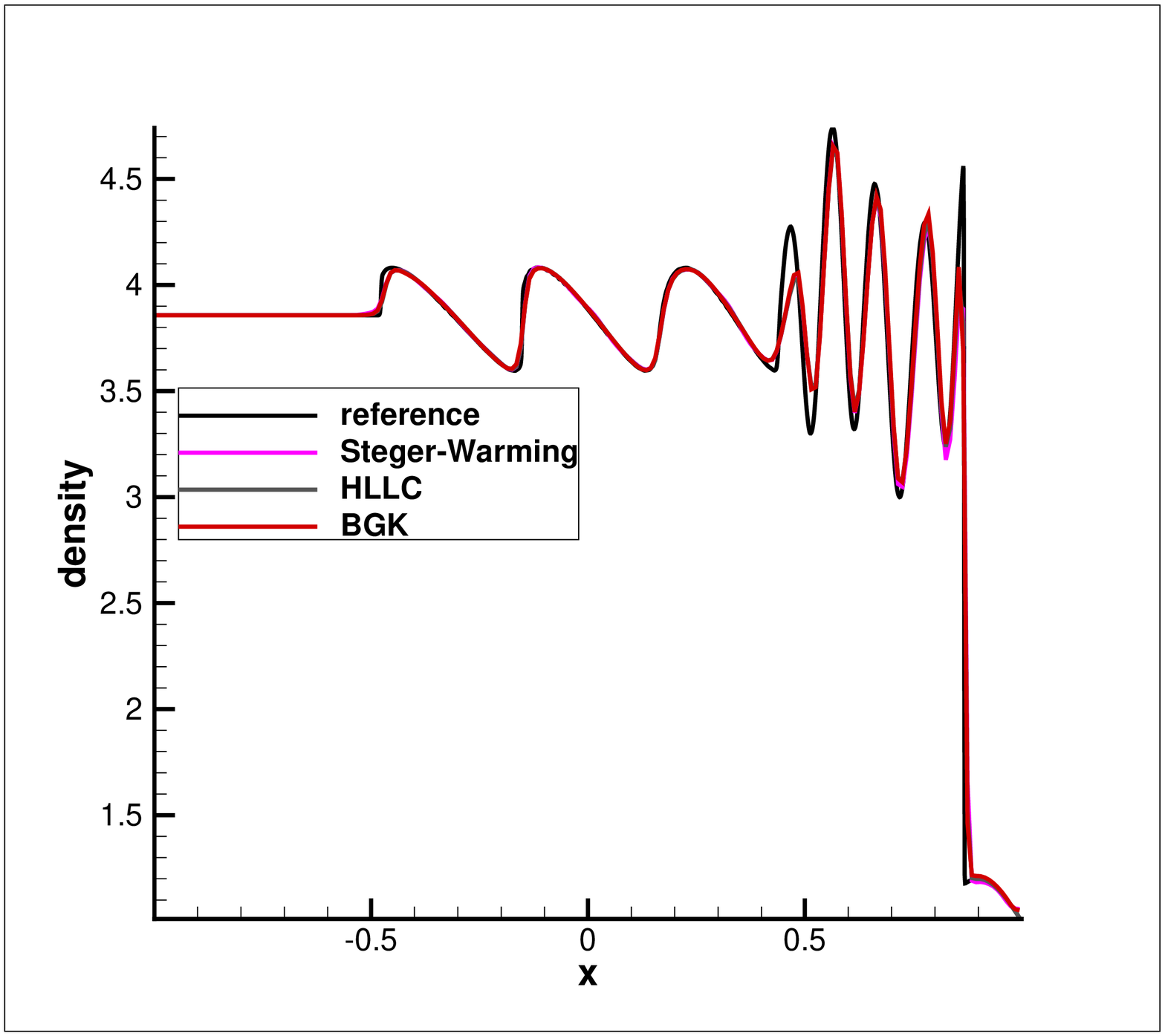}}\quad
\scalebox{0.35}{\includegraphics*[80,60][650,550]{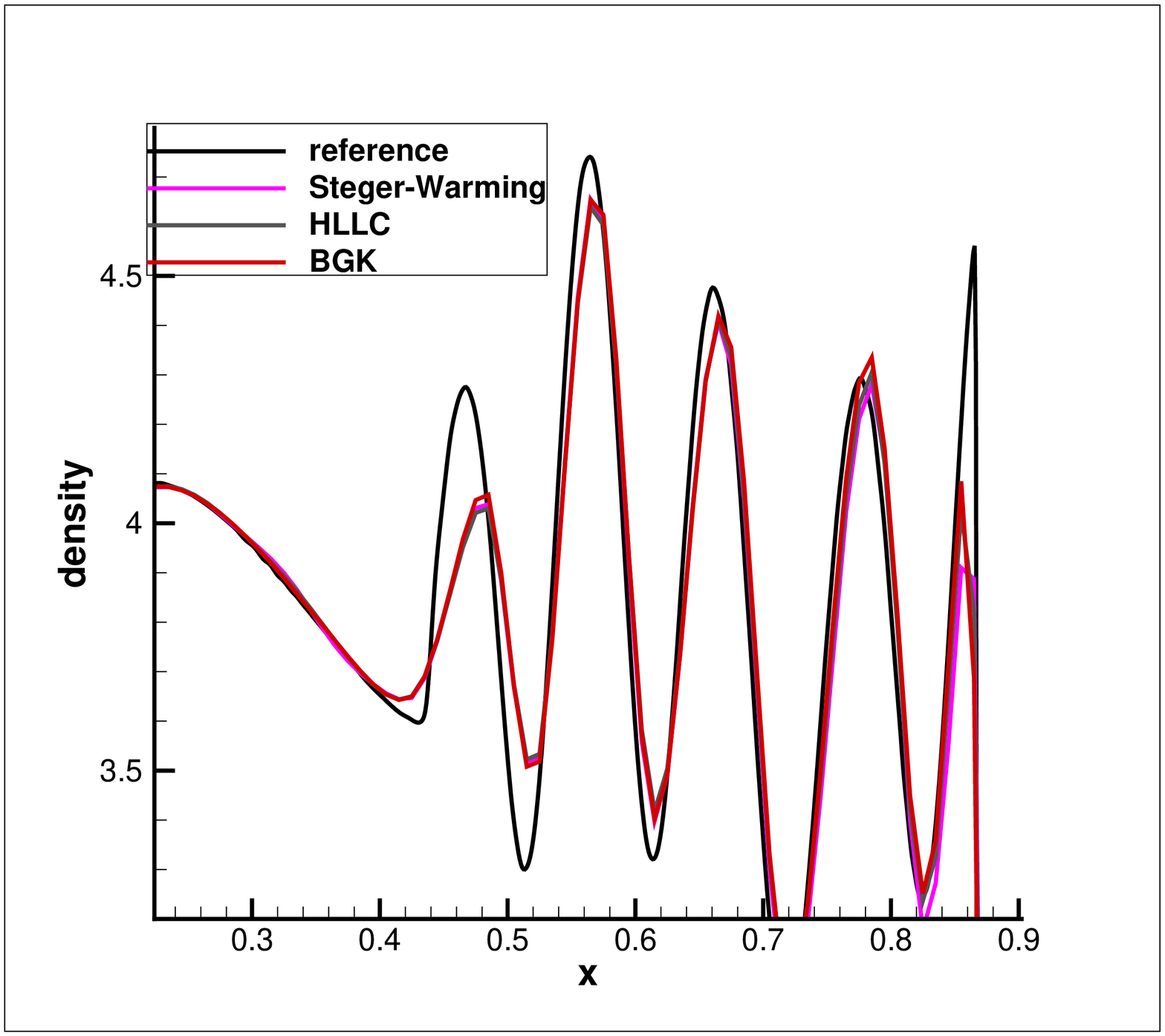}}
\end{center}
\caption{ Shu-Osher problem. density profiles for HFVS5 using
different leading terms. left : global picture, right: local
enlargement.\label{Exam4.2}}
\end{figure}

As mentioned before, it seems possible to use different leading
terms to couple with the same high derivative terms in HFVS.  We
compare three problem solver: Steger-Warming , HLLC and BGK, which
represent FVS, Godunov and GKS approaches. Fig.\ref{Exam4.1} shows
second order results and Fig.\ref{Exam4.2} shows fifth order results
of Shu-Osher problem. We can observe that all schemes using
different leading terms can work well. And there are small
differences among these scheme in low order case , while these
differences almost vanish in high order case.

These numerical results do really demonstrate the possibility: the
high order derivatives terms of HFVS can be used as a building block
to couple with any Riemann problem solvers.

\begin{figure}[!th]
\renewcommand{\baselinestretch}{0.01}
\setlength{\unitlength}{1mm} \setlength{\fboxsep}{0pt}
\begin{center}
\scalebox{0.35}{\includegraphics*[80,60][650,550]{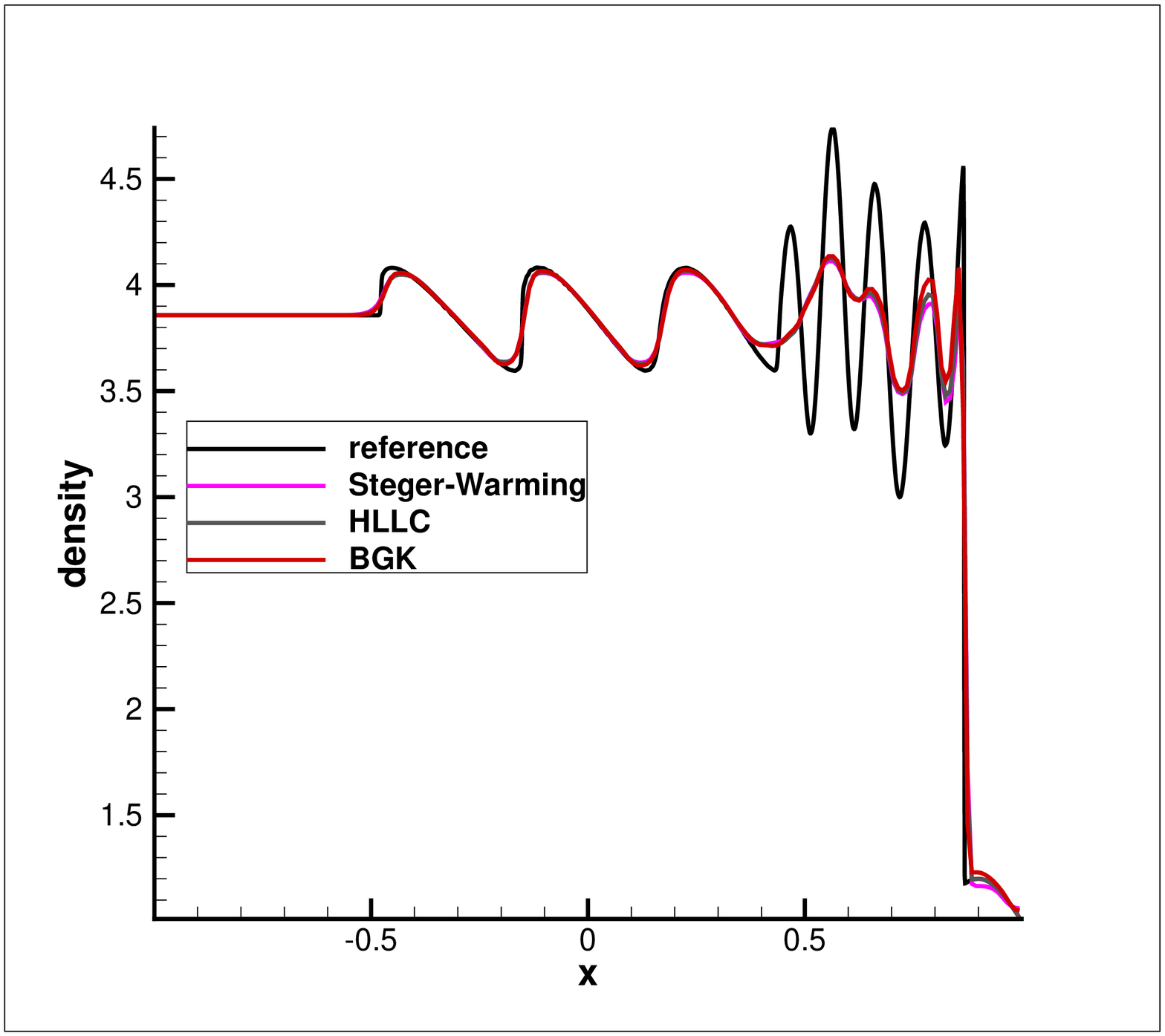}}\quad
\scalebox{0.35}{\includegraphics*[80,60][650,550]{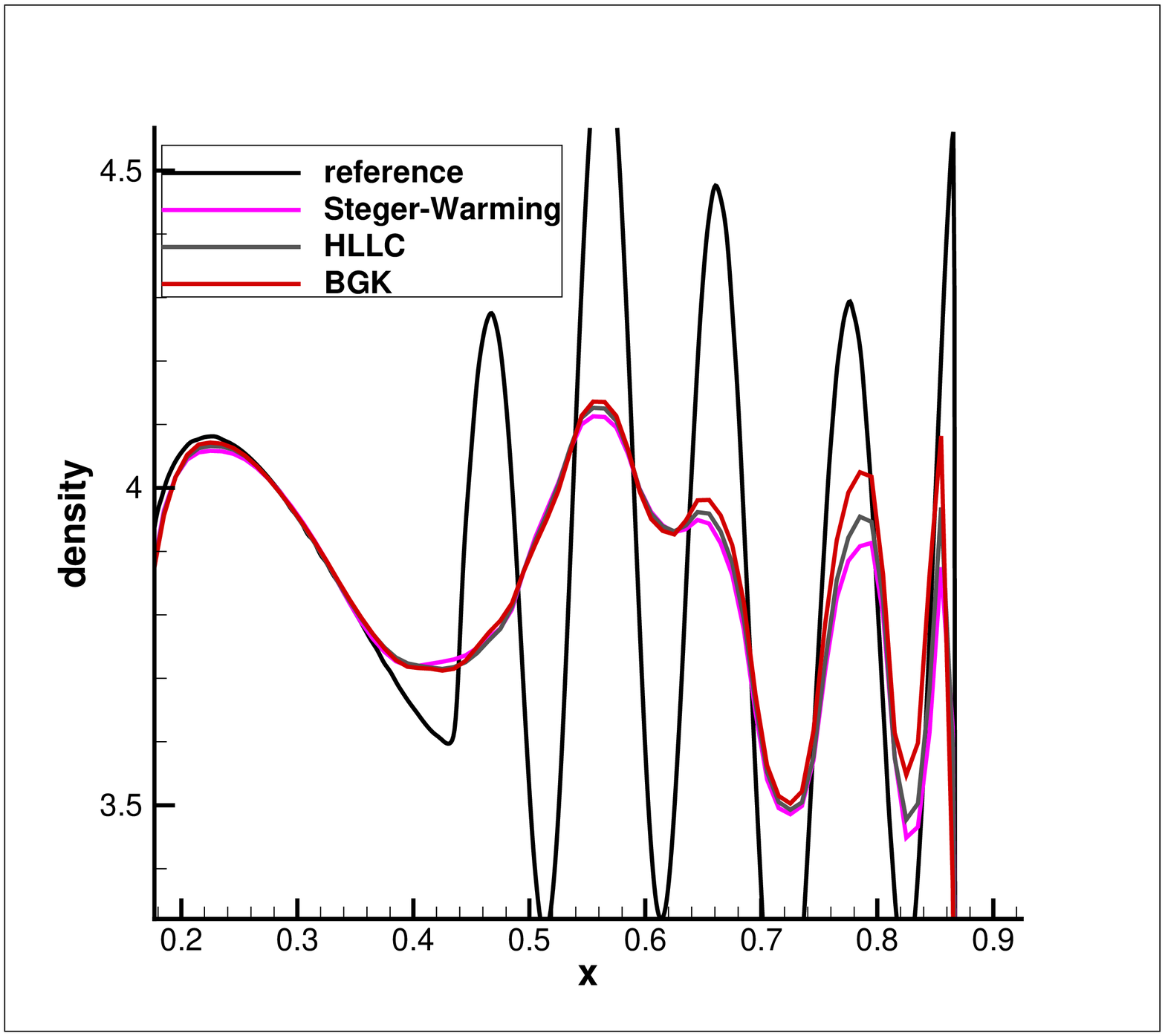}}
\end{center}
\caption{  Shu-Osher problem. density profiles for HFVS2 using
different leading terms.left : global picture, right: local
enlargement.\label{Exam4.1}}
\end{figure}

\subsubsection{Comparison between WENO and HFVS}
In fact, the detailed comparison between WENO coupling with Runge-
Kutta method with ADER can be found in \cite{Toro02}. It clearly
show the advantages of ADER scheme.

Similarly, to demonstrate the efficiency of HFVS, we present the
numerical results computed by WENO schemes coupling with Runge-Kutta
methods. We should point out that we use the same reconstruction
step in both WENO and HFVS.

Here, WENO3 and WENO5  represent third order and fifth order WENO
schemes , and RK3 is third order Runge-Kutta method.

Table 2 shows the errors and accuracy of WENO3+RK3 and WENO5+RK3
using the same computational condition as Table 1. It is clearly,
although WENO5+RK3 can get better accuracy than WENo3+RK3 , it can
also achieve third order accuracy.

Fig.\ref{Exam5.1} show the numerical results  of Shu-Osher problems,
which are computed by WENO and HFVS.  The pictures show HFVS get
better accuracy than WENO, no matter in case of third order or fifth
order. This conclusion is similar to ADER and HGKS.

Table 3 show the computational cost of WENO and HFVS. We can also
observe that HFVS almost use half CPU time of WENO. The reason is
also clear: In every computational time step, HFVS use only one
reconstruction step, while WENO should use three reconstruction
steps. Although HFVS should calculate the high order terms, it is
cheaper than reconstruction step.

\begin{table}[htbp]
 \begin{center}   \caption{Errors and accuracy for WENO}
 \label{table2}
   \begin{tabular}{|c|c|c|c|c|c|c|c|}
     \hline
     METHOD & $N$& $L_1 ERROR$ &ORDER& $L_2 ERROR$ &ORDER& $L_{\infty} ERROR$ &ORDER
     \\
     \hline
     WENO3+RK3 &  20 & 7.07E-01  &              &   7.07E-01 &              &   9.07E-01 &\\
               &  40 &  5.25E-01 &  0.429 &   5.69E-01 &  0.313 &   7.54E-01 &  0.266\\
               &  80 &  1.26E-01 &  2.058 &   1.39E-01 &  2.033 &   1.94E-01 &  1.958\\
               & 160 &  1.77E-02 &  2.831 &   1.96E-02 &  2.826 &   2.76E-02 &  2.813\\
               & 320 &  2.24E-03 &  2.982 &   2.49E-03 &  2.976 &   3.52E-03 &  2.971\\
               & 640 &  2.81E-04 &  2.994 &   3.13E-04 &  2.991 &   4.42E-04 &  2.993\\

      \hline
      WENO5+RK3& 20 &7.04E-01&        &7.04E-01&       & 7.07E-01   & \\
               & 40 &3.00E-01&  1.230 &3.40E-01&  1.050&    4.78E-01&   0.564\\
               &80  &4.42E-02&  2.762 &4.92E-02&  2.788&    6.94E-02&   2.784\\
              &160  &5.47E-03&  3.014 &6.07E-03&  3.018&    8.57E-03&   3.017\\
              &320  &6.78E-04&  3.012 &7.53E-04&  3.010&    1.06E-03&   3.015\\
              &640  &8.45E-05&  3.004 &9.39E-05&  3.003&    1.33E-04&   2.994\\
       \hline
   \end{tabular}
 \end{center}
 \end{table}

\begin{table}[htbp]
 \begin{center}   \caption{CPU time for WENO and HFVS}
 \label{table3}
   \begin{tabular}{|c|c|c|c|}
     \hline
     Method & Steps & Cpu time & Ratio \\
     \hline
     WENO3+RK3 & 1341 & 0.42    & 1 \\
     WENO5+RK3 & 1341 & 0.8268  & 1.96 \\
     HFVS3     & 1341 & 0.218   & 0.519 \\
     HFVS5     & 1341 & 0.374   & 0.89 \\
      \hline
   \end{tabular}
 \end{center}
 \end{table}

 \begin{figure}[!th]
\renewcommand{\baselinestretch}{0.01}
\setlength{\unitlength}{1mm} \setlength{\fboxsep}{0pt}
\begin{center}
\scalebox{0.5}{\includegraphics*[80,60][650,550]{pic/test0901.eps}}\quad
\scalebox{0.5}{\includegraphics*[80,60][650,550]{pic/test0902.eps}}
\end{center}
\caption{ Shu-Osher problem. Density profiles of WENO and HFVS,
left: third order schemes, right : fifth order
schemes.\label{Exam5.1}}
\end{figure}

\section{Conclusions}
In this paper, based on the idea of flux vector splitting scheme ,
we proposed a new scheme of arbitrary high order accuracy in both
space and time to solve linear and nonlinear hyperbolic conservative
laws. We start from splitting all the space and time derivatives in
the Taylor expansion of the numerical flux into two parts: one part
with positive eigenvalues, another part with negative eigenvalues.
Then all the time derivatives can be replaced by space derivatives,
according to a Lax-Wendroff procedure. A state-of-art WENO
reconstruction polynomial is used to calculate the space derivatives
. The new scheme is easy to implement, which will be very attractive
for large CFD software.

In addition, there is an interesting result: the procedure of
calculating high order terms of numerical flux can be used as a
building block to extend the current first order schemes to very
high order accuracy in both space and time.

Numerous numerical tests for linear
 and nonlinear hyperbolic conservative laws are presented to demonstrate
that new scheme is robust and can be high order accuracy in both
space and time.

%
%



\begin{thebibliography}{xx}
\bibitem{bs00} D.S.Balsara and C-W.Shu, Monotonicity preserving
weighted essentially nonoscillatory schemes with increasingly high
order of accuracy.J. Comp. Phys., 2000, {v160}, 405-452

\bibitem{grp84} M.Ben-Artzi and J.Falcovitz. . A second-order Godunov-type scheme for
compressible fluid dynamics. J. Comp. Phys. ,1984,{55}: 1¨C32.

\bibitem{Butcher}J.C.Butcher, Coefficients for the study of Runge-Kutta integration processes, J. Austral. Math. Soc. 3 (1963), 185¨C201.

\bibitem{cj09} Y.B.Chen , S.Jiang. Modified kinetic flux vector splitting schemes for compressible flows , J. Comp. Phys.
, 2009, {228(10)}:3582¨C3604

\bibitem{cj11} Y.B.Chen , S.Jiang. A non-oscillatory kinetic scheme for multi-component flows with the equation of state for a stiffened gas ,J. Comp. Phys., 2011 ,
{6}:661¨C683


\bibitem{Harten87} A.Harten,B.Engquist,S.Osher and S.R.Chakravarthy, Uniformly high oder accuracy essentially non-oscillatory schemes III, J. Comp. Phys.. 1987,{71},231-303

\bibitem{harten84}A. Harten. On a class of high resolution total-variational-stable
finite-difference schemes (with appendix by Peter D. Lax). SIAM J.
Numer. Anal., 1984 21:1{23}, .

\bibitem{weno96} G.S.Jiang and C.W.Shu . Efficient implementation of
weighed ENO schemes. J. Comput. Phys.,1996,{126}: 202¨C212.

\bibitem{hgks10} Q.Li, K.Xu, S.Fu. A high-order gas-kinetic Navier¨C Stokes
solver. J. Comput. Phys., 2010, {229}: 6715-6731

\bibitem{hgks14} N.Liu and H.Zh.Tang. A high-order accurate gas-kinetic scheme
for one- and two-dimensional flow simulation. Commun. Comput. Phys.,
2014,{15}:911-943.

\bibitem{hgks13} J.Luo, K.Xu. A high-order multidimensional gas-kinetic
scheme for hydrodynamic equations , SCIENCE CHINA Technological
Sciences, 2013, {Vol56,No.10}: 2370¨C2384

\bibitem{Roe81} P.L. Roe. Approximate riemann solver, parameter vectors and
difference schemes. J. Comput. Phys., 1981,{43}:357¨C372

\bibitem{shu-osher} C.W.Shu, S.Osher, Efficient implementation of essentially
non-oscillatory shock-capturing schemes, ii, J. Comput. Phys.,
1989,{83}:32¨C78.

\bibitem{sw81}J.L. Steger and R. Warming,Flux vector splitting of the inviscid gas
dynamic equation with application to finite differencemethods, J
Compt. Phys, 1981,{40}: 263- 293.

\bibitem{Toro02}V.A.Titarev and E.F.Toro, Arbitrary high order Godunov approach, Journal of Scientific Computing, Vol. 17. No.1-4, December 2002

\bibitem{hllc} E.F.Toro, M.Spruce and W.Speares. Restoration of the contact
surface in the Harten-Lax-van Leer Riemann solver. Shock Waves. Vol.
4, pages 25-34, 1994.

\bibitem{Toro05}E.F.Toro and V.A.Titarev, TVD fluxes for the high-order ADER schemes, Journal of Scientific Computing, Vol. 24. No.3, September 2005


\bibitem{vanLeer77} B. van Leer. Towards the ultimate conservative difference scheme.
IV: A new approach to numerical convection. J. comput. Phys.,
23:276¨C299, 1977.

\bibitem{Woodward84} Woodward and P.Colella. The numerical simulation of
two-dimensional fluid flow with strong shocks. J. Comput. Phys.,
54:115¨C173, 1984.


\end{thebibliography}
\end{document}